\documentclass[12pt,reqno]{amsart}
\usepackage{fullpage}
\usepackage{amsfonts,amsmath,amssymb,mathrsfs,amsthm,amscd} 
\usepackage{enumerate}
\usepackage{graphicx} 
\usepackage[italicdiff]{physics}
\usepackage{mathtools} 
\usepackage{array}  
\usepackage[utf8]{inputenc} 
\graphicspath{{./img/}{./pictures/}}
\renewcommand{\vec}[1]{\mathbf{#1}} 
\numberwithin{equation}{section} 
\usepackage[most]{tcolorbox} 
\newtcolorbox{mybox}[2][]
{
   colback=white,
   colframe=black,
   fonttitle=,
   coltitle=black,
   colbacktitle=white,
   enhanced,
   attach boxed title to top left={yshift=-1.2mm, xshift=2mm},
   title=#2,#1
  }

\newcommand{\eps}{\varepsilon}

\newcommand{\beq}{\begin{equation}}
\newcommand{\eeq}{\end{equation}}

\newtheorem{theorem}{\textrm{Theorem}}[section]

\newtheorem{lemma}{\textrm{Lemma}}[section]
\newtheorem{lemma1}{\textrm{Lemma}}[subsection]

\newtheorem{definition}{\textrm{Definition}}[section]
\newtheorem{definition1}{\textrm{Definition}}[subsection]

\begin{document}
	 
\begin{center}
{\large Homogenization of a coupled incompressible Stokes-Cahn-Hilliard system modeling binary fluid mixture in a porous medium}

\bigskip

\bigskip
\end{center}
\begin{minipage}[h]{0.45\columnwidth}
\begin{center}
Nitu Lakhmara\\
{\small Department of Mathematics \\IIT Kharagpur}\\
{\small WB 721302, India}\\
{\small nitulakhmara@gmail.com}
\end{center}
\end{minipage}\hspace{1cm}
\begin{minipage}[h]{0.45\columnwidth}
\begin{center}
Hari Shankar Mahato\\
{\small Department of Mathematics \\IIT Kharagpur}\\
{\small WB 721302, India}\\
{\small hsmahato@maths.iitkgp.ac.in}
\end{center}
\end{minipage}

\bigskip
\hrule
\bigskip

 	
\textbf{Abstract.} A phase-field model for two-phase immiscible, incompressible porous media flow with surface tension effects is considered. The pore-scale model consists of a strongly coupled system of Stokes-Cahn-Hilliard equations. The fluids are separated by an evolving diffuse interface of a finite width depending on the scale parameter $\varepsilon$ in the considered model. At first the well-posedness of a coupled system of partial differential equations at micro scale is investigated. We obtained the homogenized equations for the microscopic model via unfolding operator and two-scale convergence approach.
\\
\hrule
\\

\vspace{3mm}
\hrule 

\section{Introduction}
  
Fluid flow through porous media has been a subject of active research for the last four to five decades. The basic law governing the flow of fluids through porous media is Darcy's Law, which was formulated by the French civil engineer Henry Darcy in 1856 on the basis of his experiments on vertical water filtration through sand beds. A theoretical Derivation of the same can be seen in \cite{whitaker1986flow}.
In case of fluid mixtures, the dynamics of phase interface between fluids plays a central role in rheology and hydrodynamics \cite{liu2003phase,boyer2002theoretical,chaikin1995principles,feng2005energetic}.
In \cite{savin2012local} the concept of local thermodynamic equilibrium of a Gibbs interface is discussed in order to relax the global thermodynamic equilibrium assumption which is an extension towards non-equilibrium two-phase systems. 
Phase-field approach is a popular tool for the modeling and simulation of multiphase flow problems, cf. \cite{bavnas2017homogenization,kim2012phase,feng2006fully,bavnas2017numerical,anderson1998diffuse}.
The diffused interface and Cahn-Hilliard formulation of a quasi-compressible binary fluid mixture allows for topological changes of the interface, cf. \cite{lowengrub1998quasi}.
Derivation of macroscopic Cahn-Hilliard equations for perforated domains is given in \cite{schmuckupscaled}. In \cite{schmuck2013effective}, under the assumption of periodic flow and sufficiently large P$\acute{e}$clet number, a basic phase field model for a mixture of two incompressible fluids in a heterogeneous domain is upscaled via periodic homogenization. The upscaling using $\Gamma$-convergence approach for the Cahn-Hilliard equations can be found in \cite{liero2015homogenization}. 
In \cite{feng2006fully}, the author proposed fully discrete mixed finite element methods for approximating the coupled system in which, the two sets of equations are coupled through an extra phase induced stress term in the Navier–Stokes equations and a fluid induced transport term in the Cahn–Hilliard equation. In \cite{bavnas2017numerical}, a system of the Navier-Stokes equations coupled with multicomponent Cahn-Hilliard variational inequalities is considered where the existence of the weak solutions of the model is shown via a fully discrete
finite element approximation, however, the homogenization results are not shown in their work. The Cahn-Hilliard equation, in general, is given by (cf. \cite{miranville2017cahn}) 
\begin{subequations} \label{1.1}
	\begin{align} 	
	&c_t = \nabla \cdot M(c) ~ \nabla ( f(c) - \delta^2 \Delta c ), &  (t,x) \in  \mathbb{R} ^ +\times \Omega,& \label{1.1a}
	\\
	&\vec{n} \cdot \nabla c = \vec{n} \cdot M(c) \nabla ( f(c) - \delta^2 \Delta c )  = 0 , &  (t,x) \in  \mathbb{R} ^ +\times \Omega,& \label{1.1b}\\
	&c ( x , 0 ) = c_0 ( x ), &  x \in \Omega.& \label{1.1c}
	\end{align}
\end{subequations}

Here  $ 0<\delta<<1$ is a {coefficient of gradient energy}, $M = M(c)$ is a {mobility} coefficient, and $f = f (c) $ is a {homogeneous free energy}. The equation was initially developed to describe phase separation in a two component system, with $c = c(x, t)$ representing the concentration of one of the two components. Typically, the domain $\Omega$ is assumed to be a bounded domain with a {sufficiently smooth} boundary, $\partial {\Omega}$, with $\vec {n}$ in \eqref{1.1b} representing the unit exterior normal to  $\partial {\Omega}$. It is reasonable to consider evolution for times $t > 0$, or on some finite time interval $0 < t < T < \infty $.

In constant mobility quadratic polynomial case, $ M(c) = M_0 > 0 $, where $M_0$ is a constant, and $f(c) = - c + c^3 $, $f(c) = F^{\prime} (c)$, $F(c) =\frac{1}{4} (c^2 - 1)^2 $ and so Cahn-Hilliard equation reduces to 
\begin{subequations} \label{1.2}
\begin{align}
&c_t = \nabla . ~ M_0 ~ \nabla [- c + c^3 - {\delta}^2 \Delta ~ c ], &  (t,x) \in  \mathbb{R} ^ +\times \Omega , 
\label{1.2a}
\\
&\vec {n}. \nabla c = \vec {n} .  \nabla \Delta ~ c   = 0, &  (t,x) \in  \mathbb{R} ^ +\times \Omega,     
\label{1.2b}
\\
&c(x,0) = c_0 (x), & x \in \Omega.
\label{1.2c}
\end{align}
\end{subequations}

We focus on a binary-fluid model where the considered fluids are incompressible and immiscible.
The domain $\Omega \subset \mathbb{R}^n$, $ n $ = 2, 3 is occupied by the mixture. On the time interval $S = ( 0 , T )$ the model comprises a system of unsteady Stokes-Cahn-Hilliard equations
\begin{subequations} \label{1.3}
	\begin{align}
	-\mu \Delta{\vec{u}} + \nabla p &= \lambda w \nabla c & \text{ in } (0,T) \times \Omega,\label{1.3a}\\
	\nabla . \vec{u} &= 0& \text{ in }  (0,T) \times \Omega, \label{1.3b}\\
	\partial _t  c + \vec{u } . \nabla c &= \Delta w  &\text{ in } (0,T) \times \Omega, \label{1.3c}\\
	w  &= -   \Delta c  + f(c)&\text{ in } (0,T) \times \Omega, \label{1.3d}
	\end{align}
\end{subequations}
where $\vec{u}$ and $w$ are the unknown velocity and chemical potential, respectively. $\mu$ is the viscosity and $\lambda$ is the interfacial width parameter. Here $c$ represents microscopic concentration of one of the fluids with values lying in the interval $[-1, 1]$ in the considered domain and (-1, 1) within the thin diffused interface of uniform width proportional to $\lambda$. The term $f(c) = F^{\prime}(c)$, where $F$ is a homogeneous free energy functional that penalizes the deviation from the physical constraint $|c| \leq 1$. In our work we consider $F$ to be a quadratic double-well free energy functional, i.e., $
F(s) = \frac{1}{4} (s^2 - 1)^2 $. One can choose $F$ as a logarithmic or a non-smooth (obstacle) free energy functional, cf. \cite{copetti1992numerical,blowey1991cahn}.
The nonlinear term $c \nabla w$ in \eqref{1.3a} models the surface tension effects and, advection effect is modeled by the term $ \vec{u} \cdot \nabla c $ in \eqref{1.3c}. Equations \eqref{1.3a}, \eqref{1.3c}-\eqref{1.3d} represent the unsteady Stokes equations for incompressible fluid and Cahn-Hilliard equations, respectively, cf. \cite{bavnas2017homogenization}.
 
The remainder of this paper is structured as follows. In section 1.1, we  present the configuration of a periodic porous medium and the considered binary-fluid model. In section 2, notations, function space and mathematical prerequisites are gathered for the analysis which will help us to establish existence and uniqueness theorems for the pore-scale model. In  3, we will derive an anticipated macroscale model of the problem mentioned in section 1.1. In section 4, we will prove the positivity, uniqueness and existence of solution of the coupled system. Finally, in section 5, we will eventually obtain the upscaled model of the microscopic problem via two-scale convergence and unfolding operator technique. The appropriate initial and boundary conditions according to the model will be discussed further in the next section.   
\subsection{The model}
Assume that $\Omega \subset \mathbb{R}^n $, $n = 2,3$ is a bounded domain with a sufficiently smooth boundary $\partial \Omega$ and $S:=(0,T)$ denotes the time interval for any $T>0$. We consider the unit reference cell $Y \coloneqq (0,1)^n \subset  \mathbb{R}^n$. $Y_s$ and $Y_p$ represent the solid and pore part of $Y$ respectively, which are  mutually distinct, i.e., $Y_p \cap Y_s = \emptyset$, also $Y = Y_s \cup Y_p $. Solid boundary of the cell $Y$ is denoted as $\Gamma_s = \partial Y_s$, see figure 1. The domain $\Omega$ is assumed to be periodic and is covered by a finite union of the cells $Y$. In order to avoid technical difficulties, we postulate that: solid parts do not touch the boundary $\partial \Omega$, solid parts do not touch each other and solid parts do not touch the boundary of $Y$. Let $\varepsilon > 0$ be the scale parameter. We define the pore space $\Omega_p ^ \varepsilon \coloneqq \bigcup _{\mathbf{k} \in {\mathbb{Z}}^n} Y_{p_k} \cap \Omega $, the solid part as $\Omega^{\varepsilon}_s \coloneqq \bigcup _{\mathbf{k} \in {\mathbb{Z}}^n} Y_{s_k} \cap \Omega = \Omega \backslash \Omega_p^{\varepsilon}$ and $\Gamma^{\varepsilon} \coloneqq \bigcup _{\mathbf{k} \in {\mathbb{Z}}^n} \Gamma_{s_k}$, where $Y_{p_k} \coloneqq \varepsilon {Y_p + k}$, $Y_{s_k} \coloneqq \varepsilon {Y_s + k}$ and $ \Gamma_{s_k}=\bar{Y}_{p_k}\cap\bar{Y}_{s_k}$.
 Let $\chi (y)$ be the $Y$-periodic characteristic function of $Y_p$ defined by 
\begin{equation} \label{1.4}
		\chi (y) = \left\{ \begin{array}{cc}
			1 &   y \in Y^p, \\
			0 & \hspace{7mm} y \in Y-Y^p. \\
		\end{array} \right.
\end{equation}
We assume that $\Omega_p^{\varepsilon}$ is connected and has a smooth boundary. We consider the situation where the pore part $\Omega_p ^ \varepsilon$ is occupied by the mixture of two immiscible fluids separated by a microscopic interface of thickness $\eps\lambda$ ($\lambda>0$)
represented by the blue part in figure 1, and includes the effects of surface tension on the motion of the interface. 
\begin{figure}[h]
	\centering
	\includegraphics[width=13.9cm, height=5cm]{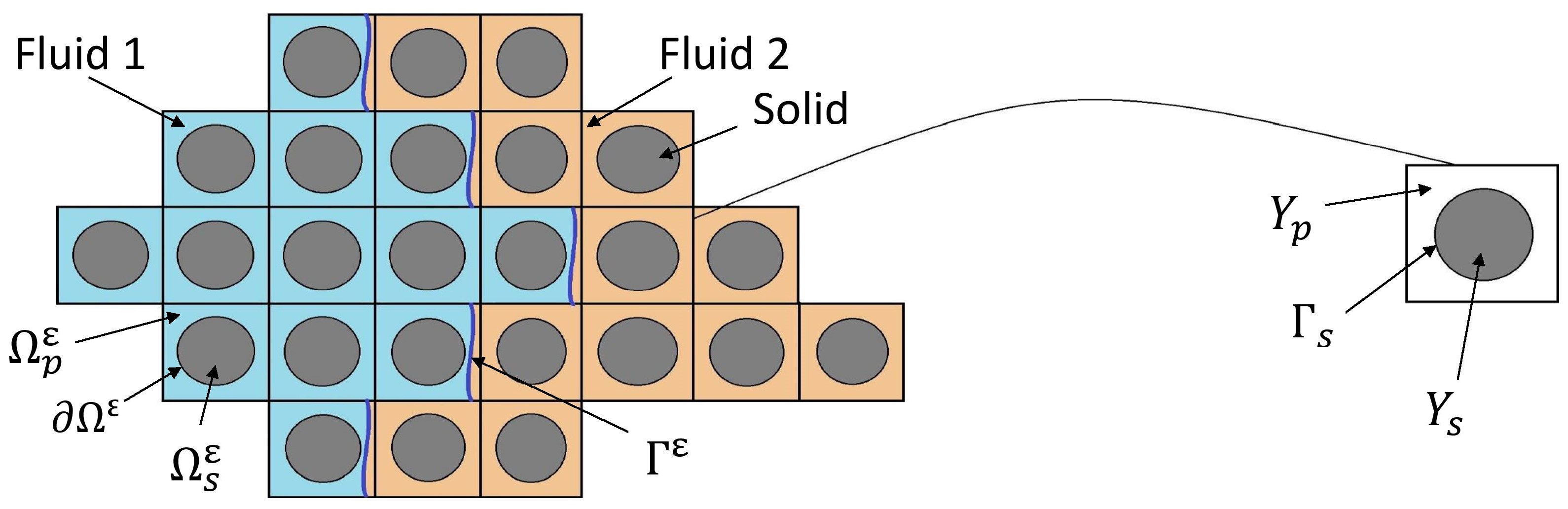}
	\caption{(left) Porous medium $\Omega = \Omega_p ^ \epsilon\cup {\Omega}_s ^{\epsilon}$ as a periodic covering of the reference cell $Y = Y_p \cup Y_s$ (right). The blue interface is the microscopic boundary between two fluid phases occupying the the pore space $\Omega_p^ \epsilon$.}	
\end{figure} 
We model the flow of the fluid mixture on the pore-scale using a phase-field approach motivated by the Stokes-Cahn-Hilliard system \eqref{1.3}. The velocity of the fluid mixture is assumed to be $\mathbf{u}^\varepsilon = \mathbf{u}^\varepsilon (t,x)$,  $(t,x) \in S \times \Omega_p^ \varepsilon $ which satisfies the stationary Stokes equation. The order parameter $c^ \varepsilon$ plays the role of microscopic concentration and the chemical potential $w^ \varepsilon$ satisfy the Cahn-Hilliard equation. $p^ \varepsilon$ is the fluid pressure. The term $\eps\lambda c^ \varepsilon \nabla w^ \varepsilon$ models the surface tension forces which acts on the microscopic interface between the fluids. Fluid density is taken to be $1$. Then, the  Stokes-Cahn-Hilliard system of equations is given by 
  
\begin{mybox}{Pore scale equations}
\begin{subequations}  \label{1.5}
	\begin{align}
			- \mu \varepsilon^2 \Delta \vec{u^\varepsilon} + \nabla p^\varepsilon &  = - \varepsilon  \lambda c^\varepsilon \nabla w^\varepsilon &   S \times \Omega_p^ \varepsilon,  \label{1.5a}\\ 
			\nabla . \vec{u^\varepsilon} & = 0 &  S \times \Omega_p^ \varepsilon,  \label{1.5b}\\
			\vec{u^\varepsilon} & = 0 &  S \times \partial\Omega_p^ \varepsilon,   \label{1.5c}\\
			\partial _t  c^\varepsilon + \varepsilon ^ \beta \vec{u^\varepsilon} . \nabla c^\varepsilon &  = \varepsilon ^ \alpha \Delta w^\varepsilon &  S \times \Omega _p^ \epsilon, \label{1.5d}\\
			w^\varepsilon & = - \varepsilon ^ \gamma \Delta c^\varepsilon + f(c^\varepsilon) &  S \times \Omega_p^ \epsilon,   \label{1.5e}\\
		\partial_n  c^\varepsilon& = 0 &   S \times \partial\Omega_p^ \varepsilon,  \label{1.5f}\\  
		\partial_n  w^\varepsilon& = 0 &   S \times \partial\Omega_p^ \varepsilon,  \label{1.5g}\\
			c^\varepsilon (0,x) & = c_0(x) &  \Omega_p^{\varepsilon}, \label{1.5h}
		\end{align}
\end{subequations}
\end{mybox}
where $\frac{\partial   c^\varepsilon}{\partial \vec{n}}=\partial_n c^\varepsilon$ and $f(s) = s^3 - s = F^{\prime}(s) = \frac{1}{4} (s^2 - 1)^2 $ is the double-well free energy. The above scaling for the viscosity is such that the velocity $\vec{u^\varepsilon}$ has a nontrivial limit as $\varepsilon$ goes to zero. Also, $ 0 \leq \alpha , \beta, \gamma \leq 2 $ where $ \alpha, \beta, \gamma \in \mathbb{R}$. We denote (1.5a)-(1.5h) by $\mathcal{(P^\eps)}$.

	
\section{Notations and function space setting}
Let $\theta \in [0,1]$ and $1\le r,s\le \infty$ be such that $\frac{1}{r}+\frac{1}{s}=1$. Assume that $\Xi\in \{\Omega, \Omega_{p}^{\varepsilon}, \Omega_{s}^{\varepsilon}\}$ and $l\in\mathbb{N}_0$, then as usual $L^{r}(\Xi)$ and $H^{l,r}(\Xi)$ denote the Lebesgue and Sobolev spaces with their usual norms and they are denoted by $||.||_r$ and $||.||_{l,r}$. Similarly, $C^{\theta}(\bar{\Xi})$, $(\cdot,\cdot)_{\theta,r}$ and $[\cdot,\cdot]_{\theta}$ are the H\"older, real- and complex-interpolation spaces respectively endowed with their standard norms, for definition confer \cite{Eva98, Lun95}. $C^{\alpha}_{\#}(Y)$ denotes the set of all \textit{Y-periodic $\alpha$-times} continuously differentiable functions in $y$ for $\alpha\in \mathbb{N}$. In particular, $C_{\#}(Y)$ is the space of all the \textit{Y-periodic} continuous function in $y$. The $C^\infty$-spaces are as usual equipped with their \textit{maximum norm} whereas the space of all continuous functions $C(\Xi)$ is furnished with \textit{supremum norm}, cf. in \cite{Eva98}.
The extension and restriction operators are denoted by $E$ and $R$, respectively. The symbol $(.,.)_{H}$ represents the \textit{inner product} on a \textit{Hilbert space H} and $||.||_{H}$ denotes the corresponding norm. For a Banach space $X$, $X^{*}$ denotes its dual and the duality pairing is denoted by $\langle.\; ,\; .\rangle_{X^{*}\times X}$. By classical trace theorem on \textit{Sobolev space} $H^{1,2}_{0}(\Xi)^*=H^{-1,2}(\Xi)$. The symbols $\hookrightarrow$, $\hookrightarrow\hookrightarrow$ and $\underset{\hookrightarrow}{d}$ denote the continuous, compact and dense embeddings respectively. By definition, the surface area of $\Gamma^{\varepsilon}$ increases proportionally to $\frac{1}{\varepsilon}$, i.e. $|\Gamma^{\varepsilon}|\to \infty$ as $\varepsilon\to 0$. 
We define the function spaces: 
	\begin{center}
		
		$\mathbf{H}^1 (\Omega) =  H^1 (\Omega)^n ,\quad \mathbf{H}^1_0 (\Omega) =  H^1_0 (\Omega)^n ,\quad \mathfrak{U}^\varepsilon :=\mathbf{H}^1_{div} (\Omega) = \{\eta : \eta \in \mathbf{H}^1_0 (\Omega), \nabla \cdot \eta = 0 \}$,
		\vspace{1mm}\\
		$\mathfrak{C}^\varepsilon = \{ c^ \varepsilon : c^ \varepsilon \in L^{\infty}(S ; H^1(\Omega_p ^ \varepsilon)), \partial_t c^ \varepsilon \in L^2(S ; H^1(\Omega_p ^ \varepsilon)^*) \}$,
	    \vspace{1mm}\\ 
	    $\mathfrak{W}^\varepsilon = L^2( S; H^1(\Omega_p ^ \varepsilon))$\text{ and }$L^2_0(\Omega)=\{\phi\in L^2(\Omega): \int_{\Omega}\phi\,dx=0.\}$.
	\end{center}
	
We choose $ \vec{u^\varepsilon} \in  \mathfrak{U}^\varepsilon  $, $c^ \varepsilon \in \mathfrak{C}^\varepsilon$, $w^ \varepsilon \in  \mathfrak{W}^\varepsilon$ and  $p^\varepsilon \in L^2(S\times \Omega_p ^ \varepsilon)$. We will now state few results and lemmas which are used in this paper and proofs of these can be found in literature.  
 
\begin{lemma}[cf. theorem II.1.1 in \cite{Sho97}]
		Let $E$ be a Banach space and $E_0$ and $E_1$ be reflexive spaces
		with $E_0 \subset E \subset E_1$. Suppose further that $E_0 \hookrightarrow \hookrightarrow E \hookrightarrow E_1$. For $1 < p, q < \infty$ and $0 < T < 1$ define $ X \coloneqq \{ u \in L ^ p ( S ; E_0 ) : \partial_t u \in L^q ( S ; E_1 ) \} $. Then $ X \hookrightarrow \hookrightarrow L ^p ( S ; E )  $.
\end{lemma}
 
	
\begin{lemma}[Restriction theorem, cf. \cite{Rad92,All92,allaire1993school}]
	There exists a linear restriction operator $ R^{\varepsilon} : L^2(S;H_0^1 ( \Omega ))^d \longrightarrow L^2(S;H_0^1 ( \Omega ^{ \varepsilon }_p ))^d $ such that $ R^{\varepsilon} u(x) = u(x) | _{\Omega ^{\varepsilon}_p} $ for $ u \in L^2(S;H_0^1 ( \Omega ))^d $ and $ \nabla \cdot R^{\varepsilon} u = 0 $ if $ \nabla \cdot R^{\varepsilon} u = 0 $ if $ \nabla \cdot u = 0 $. Furthermore, the restriction satisfies the following bound 		
	\begin{equation*}
			\norm{  {R}^{\varepsilon} u }_{L^2(S\times\Omega ^{\varepsilon}_p)} + \varepsilon \norm{ \nabla {R}^{\varepsilon} u }_{L^2(S\times\Omega ^{\varepsilon}_p)} \leq  C ( \norm{ u }_{L^2 (S\times \Omega )} + \varepsilon \norm{ \nabla  {u} }_{ L^2(S\times\Omega) } ),
	\end{equation*}		
where $C$ is independent of $\varepsilon$.
\end{lemma}	
Similarly, one can define the extension operator from $S\times \Omega ^{ \varepsilon }_p$ to $S\times\Omega$, cf. \cite{Rad92,Pet07,mahato2013homogenization,All92,allaire1993school}.

		


		
\setcounter{definition}{4}
\begin{definition}[Two-scale convergence, cf. \cite{Ngu89,All92,Rad92}]
	A sequence of functions $(u^{\varepsilon})_{\varepsilon > 0}$ in $ L^{p} ( S \times \Omega ) $ is said to be two-scale convergent to a limit $ u \in L^{p} ( S \times \Omega \times Y ) $ if
	\begin{equation*}
			\lim_{\epsilon \rightarrow 0} \int_{ S \times \Omega } u^{\varepsilon} (t , x) \phi (t, x, \frac{x}{\varepsilon}) \,dx \,dt = \int_{ S \times \Omega \times Y} u(t, x, y) \phi(t, x, y) \,dx \,dt \,dy
	\end{equation*}
	for all $ \phi \in L^q( S \times \Omega ; C_{\#}( Y ) ) $.
\end{definition}
	
By $ \overset{2}{\rightharpoonup} $, $ \overset{w}{\rightharpoonup} $ and $ \rightarrow $ we denote the two-scale, weak and strong convergence of a sequence respectively. The proof of next lemma is given in \cite{Ngu89,All92,Rad92}.

\begin{lemma}
For $\eps>0$, let $(u^{\varepsilon})_ {\varepsilon > 0}$ be a sequence of functions, then the following holds:
\begin{itemize}
\item[(i)] for every bounded sequence $(u^{\varepsilon})_ {\varepsilon > 0}$ in $ L^p ( S \times \Omega ) $ there exists a subsequence $(u^{\varepsilon})_ {\varepsilon > 0}$ (still denoted by same symbol) and an $ u \in L^p ( S \times \Omega \times Y) $ such that $  u^{\varepsilon} \overset{2}{\rightharpoonup} u  $.
\item[(ii)] let $u^{\varepsilon} \to u $ in $ L^p ( S \times \Omega ) $, then $  u^{\varepsilon} \overset{2}{\rightharpoonup} u  $.
\item[(iii)] let $(u^{\varepsilon}) _ {\varepsilon > 0}$ be a sequence in $ L^{p}(S; H^{1,p}(\Omega) ) $ such that $  u^{\varepsilon} \overset{w}{\rightharpoonup} u  $ in $ L^{p}(S; H^{1,p}(\Omega) ) $. Then $ u^{\varepsilon} \overset{2}{\rightharpoonup} u $ and there exists a subsequence $u^{\varepsilon} _ {\varepsilon > 0}$, still denoted by same symbol, and an $ u_1 \in L^p (S\times \Omega ; H_{\#}^{1,p} (Y) ) $ such that $ \nabla_{x} u^{\varepsilon} \overset{2}{\rightharpoonup} \nabla_{x} u + \nabla_{y} u_1 $.
\item[(iv)] let $(u^{\varepsilon}) _ {\varepsilon > 0}$ be a bounded sequence of functions in $L^p (S \times \Omega )$ such that $ \varepsilon \nabla u^{\varepsilon} $ is bounded in $L^{p} (S \times \Omega)^n$. Then there exist a function $ u \in L^p (S \times \Omega ; H_{\#}^{1,p} (Y) ) $ such that $ u^{\varepsilon} \overset{2}{\rightharpoonup} u $, $ \varepsilon \nabla_{x} u^{\varepsilon} \overset{2}{\rightharpoonup} \nabla_{y} u $.
\end{itemize}
\end{lemma}
	
	
	

\subsection{Periodic Unfodling}
Arbogast, Douglas and Hornung in \cite{ADH90} introduced the concept of \textit{dilation operator} to study the homogenization on periodic domains with double porosity. This method is further used in \cite{BLM96}, \cite{RJ07}, \cite{CP08} and references therein. Consequently, the idea of dilation operator is extended by Cioranescu, Damlamian and Griso (cf. \cite{CDG02}, \cite{CDG08}) to examine the homogenization problems on periodic domains under the name of \textit{periodic unfolding}. The unfolding operator technique has proven to be useful when we deal with the homogenization problems with nonlinear terms. We start by the notations and terminologies required to define the unfolding operator.
Let $\Omega$, $\Omega_p^\varepsilon$, $\Omega_s^\varepsilon$ and $\Gamma^\varepsilon$ be defined as in section 1.1. For any $z\in \mathbb{R}^{n}$, suppose $[z]$ denotes the unique integer combination $\sum_{j=i}^{n}k_{j}e_{j}$ of $e_{j}$ such that $z-[z]$ lies in $Y$ (see figure 2 in \cite{CDZ06}) and we set $\left\{z\right\}=z-[z]\quad \textnormal{ for a.e. }z\in \mathbb{R}^{n}.$ Thus for any $x\in \mathbb{R}^{n}$ and $\varepsilon >0$, we have $x=\varepsilon\left(\left[\frac{x}{\varepsilon}\right]+\left\{\frac{x}{\varepsilon}\right\}\right)\quad \textnormal{ a.e. }x\in \mathbb{R}^{n}.$ Setting $\Theta^\varepsilon=\left\{\xi\in \mathbb{Z}^{n}:\varepsilon(\xi+Y)\subset \Omega\right\},\; \hat{\Omega}^\varepsilon=\textnormal{ int}\left\{\cup_{\xi\in \Theta^\varepsilon}\varepsilon (\xi+\bar{Y})\right\} \text{ and }\Upsilon^\varepsilon=\Omega-\hat{\Omega}^\varepsilon.$



\begin{definition1}\label{def2.3.1}
Assume that $1\le r\le \infty$. Let $u^\varepsilon \in L^{r}( S \times \Omega )$ such that for every $t$, $u^\varepsilon(t)$ is extended by zero outside of $\Omega$. We define the unfolding operator $T^{\varepsilon}:L^{r}( S \times \Omega)\to L^{r}( S \times \Omega \times Y)$ as
\begin{subequations}
\begin{align}
T^{\varepsilon}u^\varepsilon(t, x,y)&=u^\varepsilon\left(t,\varepsilon \left[\frac{x}{\varepsilon}\right]+\varepsilon y\right)&& \textnormal{for a.e. }(t, x, y)\in S \times {\Omega} \times Y,\label{2.30a}\\
  &=0&&\textnormal{otherwise}.\label{2.30b}
\end{align}\end{subequations}
\end{definition1}
Based on definition 2.3.1 the unfolding operator $T^\varepsilon$ preserves the integral and the norms on the domain $\Omega_p^\varepsilon$
\[\int_{\Omega}u^\varepsilon(x)\,dx=\frac{1}{|Y|}\int_{\Omega\times Y}T^\varepsilon u^\varepsilon(x,y)\,dx\,dy\text{ and }||u^\varepsilon||_{L^r(\Omega )}=\frac{1}{|Y|^\frac{1}{r}}||T^\varepsilon u^\varepsilon||_{L^r(\Omega\times Y)}.\]
$T^\varepsilon $ is linear. We also note that $\nabla_y T^\varepsilon u^\varepsilon(x,y) = \varepsilon T^\varepsilon(\nabla_x u^\varepsilon)(x,y)$ and  $\Delta_y T^\varepsilon u^\varepsilon(x,y) = \varepsilon^2 T^\varepsilon$ $(\Delta_x  u^\varepsilon)(x,y)$. Based on definition 2.3.2, the following properties of $T^\varepsilon$ can be proved (cf. \cite{CDG02}, \cite{CDG08} and \cite{CDZ06}):
\setcounter{lemma1}{2}
\begin{lemma1}\label{thm3.6.4}
Let $1<p<\infty$, then the operator $T^\varepsilon$ has the following properties:
\begin{enumerate}
\item[(i)] If $u^\varepsilon\in L^{p}(S\times \Omega)$, then  $T^\varepsilon u^\varepsilon(t,x, \left\{\frac{x}{\varepsilon}\right\})=u^\varepsilon(t,x)$,\textnormal{ for every }$t\in S$ and $x\in \Omega$.
\item[(ii)] Let $u^\varepsilon,v^\varepsilon \in L^{p}(S\times \Omega)$, then $T^\varepsilon(u^\varepsilon v^\varepsilon)=T^{\varepsilon}(u^\varepsilon)T^\varepsilon(v^\varepsilon)$.
\item[(iii)] Let $u^\varepsilon\in L^{1}(S\times \Omega)$, then ${\int}_{0}^{T}{\int}_{\Omega}u^\varepsilon(t,x)\,dx\,dt=\frac{1}{\varepsilon|Y|}{\int}_{0}^{T}{\int}_{\Omega \times Y}T^\varepsilon(u^\varepsilon(t,x,y))\,dx\,d y\,dt$.
\item[(iv)] Let $u^{\varepsilon}\in L^{p}(S\times \Omega)$, then ${\int}_{0}^{T}{\int}_{\Omega}{\int}_{Y}\left|T^\varepsilon u^\varepsilon(t,x,y)\right|^{p}\,dx\,dy\,dt=\varepsilon|Y|{\int}_{0}^{T}{\int}_{\Omega}\left|u^\varepsilon(t,x)\right|^{p}\,dx\,dt$.
\item[(v)] $u^\eps \overset{2}{\rightharpoonup}u$ if and only if $T^\eps u^\eps \overset{w}{\rightharpoonup}u$.
\end{enumerate}
\end{lemma1}

For the following definitions and results, interested reader can refer to \cite{francu2007two,francuu2012some,francuu2010modification} and references therein.  
\setcounter{definition1}{3}
\begin{definition1}[cf. definition 4.5 in \cite{francuu2012some}]\label{def2.3.1}
Assume that $1\le r\le \infty, u^\varepsilon \in L^{r}(S\times \Omega)$ and $T^\varepsilon $ is defined as in definition \ref{def2.3.1}. Then we say that:\\
(i) $u^\varepsilon$ is weakly two-scale convergent to a limit $u_0\in L^r(S\times\Omega\times Y)$ if $T^\varepsilon u^\varepsilon$ converges weakly to $ u_0$ in $L^r(S\times\Omega\times Y)$.\\
(ii) $u^\varepsilon$ is strongly two-scale convergent to a limit $u_0\in L^r(S\times\Omega\times Y)$ if $T^\varepsilon u^\varepsilon$ converges strongly to $ u_0$ in $L^r(S\times\Omega\times Y)$.
\end{definition1}

\setcounter{lemma1}{4}
\begin{lemma1}\label{thm3.6.2}
Let $\left(u^{\varepsilon} \right)_{\varepsilon >0}$ be a bounded sequence in $L^{r}(S\times \Omega)$. Then the following statements hold:
\begin{enumerate}
\item[(a)] if $u^\eps\overset{2}{\rightharpoonup}u$, then $T^{\varepsilon} u^{\varepsilon}\overset{w}{\rightharpoonup} u$, i.e. $u^\varepsilon$ is weakly two-scale convergent to a $u$.
\item[(b)] if $u^\eps\to u$ , then $T^{\varepsilon} u^{\varepsilon}\to u$, i.e. $u^\varepsilon$ is strongly two-scale convergent to  $u$.
\end{enumerate}
\end{lemma1}
\begin{proof}
(a) Since $u^\eps\overset{2}{\rightharpoonup}u$, then by definition we have $\underset{\varepsilon \to 0}{\lim}\;\;\;\int_{0}^{T}\int_{\Omega}u^{\varepsilon}(t,x)\phi(t,x,\frac{x}{\varepsilon})\,dx\,dt=\int_{0}^{T}\int_{\Omega}\int_{Y}u(t,x,y)\phi(t,x,y)\,dx\,d y\,dt $ for all $\phi \in C( S \times \bar{\Omega} ; C_{\#}(Y))$. Now, from the linearity and integral preserving property of $T^\eps $ (note that $|Y|=1$), we have $\int_{S\times\Omega}u^\eps\phi\,dx \,dt=\int_{S\times\Omega\times Y}T^\eps (u^\eps\phi)\,dx\,dy\,dt=  \int_{S\times\Omega\times Y}T^\eps u^\eps T^\eps \phi\,dx\,dy\,dt=\int_{S\times\Omega\times Y}T^\eps  u^\eps \phi\,dx\,d y\,dt .$ This gives $\lim_{\eps\to 0}\int_{S\times\Omega\times Y}T^\eps  u^\eps \phi\,dx\,d y\,dt =\lim_{\eps\to 0}\int_{S\times \Omega}u^\eps\phi\,d  x\,dt=\int_{S\times\Omega\times Y}u \phi\,dx\,d y\,dt$, i.e. $T^\eps u^\eps $ is weakly convergent to $u$, i.e. $u^\eps $ is weakly two-scale convergent to $u$.\\
(b) We note that by integral preserving property, $||T^\eps u^\eps- u||_{L^2(S\times\Omega\times Y)}^2=(\eps|Y|)^{\frac{1}{2}}||u^\eps- u||_{L^2(S\times\Omega)}^2\to 0\text{ as }\eps \to 0. $ This implies $T^\eps u^\eps $ is strongly convergent to $u$, i.e. $u^\eps $ is strongly two-scale convergent to $u$.
\end{proof}
\begin{lemma1}[cf. theorem 6.2 in \cite{francuu2012some}]\label{thm3.6.2}
Let $\left(u^{\varepsilon} \right)_{\varepsilon >0}$ be strongly two-scale convergent to $u_0$ in $L^{r}(S\times \Omega\times \Gamma)$ and $\left(v^{\varepsilon} \right)_{\varepsilon >0}$ be weakly two-scale convergent to $v_0$ in $L^{s}(S\times \Omega\times \Gamma)$. If the exponents $r,s,\nu\ge 1$ satisfy $\frac{1}{r}+\frac{1}{s}=\frac{1}{\nu}$, then the product $(u^\varepsilon v^\varepsilon)_{\varepsilon>0}$ two-scale converges to the limit $u_0v_0$ in $L^\nu(S\times\Omega\times Y)$.
In particular, for any $\phi\in L^\mu(S\times\Omega)$ with $\mu\in (1,\infty)$ such that $\frac{1}{\nu}+\frac{1}{\mu}=1$ we have
\[
\int_{S\times\Omega}u^\varepsilon(t,x) v^\varepsilon(t,x) \phi(t,x)\,dx\,dt\overset{\varepsilon\to 0}{\longrightarrow}\int_{S\times \Omega\times  Y} u_0(t,x,y)v_0(t,x,y)\phi(t,x)\,dx\,d y\,dt.
\]
\end{lemma1}

Before we proceed with the weak formulation, we make following assumptions for the sake of analysis of $\mathcal{(P^\eps)}$.

\begin{enumerate}
\item[\textbf{A1.}] for all $x\in \Omega$, $ \vec{u_0}$, $c_0$ and $w_0\ge 0$.
\item[\textbf{A2.}] $\vec{u_0}\in L^\infty(\Omega)\cap H^1(\Omega)$, $c_0\in L^\infty(\Omega)\cap H^1(\Omega)$ and $w^0\in L^\infty(\Omega)\cap H^1(\Omega)$ such that $\sup_{\eps>0}||\vec{u_0}||_{L^\infty(\Omega)\cap H^1(\Omega)}<\infty,$ $ \sup_{\eps>0}||c_0||_{L^\infty(\Omega)\cap H^1(\Omega)}<\infty,$ $ \sup_{\eps>0}||{w_0}||_{L^\infty(\Omega)\cap H^1(\Omega)}<\infty$.
\item[\textbf{A3.}] $p^\eps \in L^2(S;H^1(\Omega_p^\eps))$ such that $\sup_{\eps>0}||p^\eps||_{L^2(S;H^1(\Omega_p^\eps))}<\infty$.
\item[\textbf{A4.}] $\alpha = 2$, $ \beta = 1 $ and $ \gamma = 0 $.
\end{enumerate}


\subsection*{2.10. Weak formulation of $\mathcal{(P^\eps)}$} Let the assumptions A1 - A4 be satisfied. A triple $(\vec{u^\varepsilon}$, $c^ \varepsilon$,  $w^ \varepsilon) \in \mathfrak{U}^\eps\times \mathfrak{C}^\eps\times \mathfrak{W}^\eps$ is said to be the weak solution of the model $\mathcal{(P^\eps)}$ such that $ (\vec{u^\varepsilon}, c^\varepsilon, w^\varepsilon)(0,x) = (\vec{u}_0, c_0, w_0)(x)$ for all $x\in  \Omega $, and
\begin{subequations} \label{2.1}
\begin{align} 
&\mu \varepsilon ^2 \int_{S \times \Omega_p^ \varepsilon} \nabla \vec{u^\varepsilon} : \nabla \eta \,dx \,dt = - \varepsilon  \lambda \int_{S \times \Omega_p^ \varepsilon} c^{\varepsilon} \nabla w^{\varepsilon} \cdot \eta \,dx \,dt , 
\label{2.1a}\\
&\int_{S} \langle \partial_t c^{\varepsilon} , \phi \rangle \,dt -  \varepsilon   \int_{S \times \Omega_p^ \varepsilon} c^{\varepsilon} \vec{u^{\varepsilon}} \cdot \nabla \phi \,dx \,dt + \varepsilon^2   \int_{S \times \Omega_p^ \varepsilon} \nabla
w^{\varepsilon} \cdot \nabla \phi \,dx \,dt = 0, 
\label{2.1b}\\
&\int_{S \times \Omega_p^ \varepsilon} w^{\varepsilon} \psi \,dx \,dt =  \int_{S \times \Omega_p^ \varepsilon}  \nabla c^{\varepsilon} \cdot \nabla \psi \,dx \,dt +  \int_{S} \langle  f(c^{\varepsilon}) , \psi \rangle \,dx \,dt,
\label{2.1c}
\end{align}
\end{subequations}
for all $ \eta \in L^2(S;\mathbf{H}^1_{div} (\Omega_p ^ \varepsilon)) $ and $ \phi $, $ \psi  \in  L^2(S; {H}^1 (\Omega_p ^ \varepsilon)) $.


We are now going to state the two main theorems of this paper which are given below.

\begin{theorem}
Let the assumptions A1 - A4 be satisfied, then there exists a unique positive weak solution $(\vec{u^\varepsilon}$, $c^ \varepsilon$,  $w^ \varepsilon) \in \mathfrak{U}^\eps \times \mathfrak{C}^\eps\times \mathfrak{W}^\eps$ of the problem $\mathcal{(P^\eps)}$ which satisfies
\begin{align} \label{2.2}
&\sqrt{\mu} \varepsilon   \norm{\nabla \vec{u^\varepsilon} }_{L^2(S \times {\Omega}_p^{\varepsilon})}   +  \sqrt{\varepsilon \lambda}  \norm{ \nabla w^{\varepsilon} }_{L^2(S\times {\Omega}_p^{\varepsilon})} + \sqrt{\frac{\lambda}{2}} \norm{ \nabla c^{\varepsilon} }_{L^{\infty} (S); L^2 ({\Omega}_p^{\varepsilon}) )}+\norm{ w^{\varepsilon} }_{L^2(S \times {\Omega}_p^{\varepsilon})}
\notag\\
&\qquad\quad+\norm{ c^{\varepsilon} }_{L^{\infty}(S ; L^4 ({\Omega}_p^{\varepsilon}))}+\norm{\partial_t c^{\varepsilon}}_{L^2(S; H^1 ({\Omega}_p ^ \varepsilon)^* )}+\norm{\vec {u^{\varepsilon}}}_{L^4(\Omega_p^{\varepsilon})}\le C<\infty  \quad \forall \varepsilon,
\end{align}
where the constant $C$ is independent of $\eps$.
\end{theorem}


\begin{theorem}[Upscaled Problem $\mathcal{(P)}$]\label{thm2.2.2}
There exists $(\vec{u} ,c ,w) \in \mathfrak{U}  \times \mathfrak{C} \times \mathfrak{W} $ which satisfies 
\begin{subequations} \label{2.3}
\begin{align}
  - \mu \Delta_y \vec{u} ( x,y )  + \nabla_y p_1 (x,y)  + \nabla_x p ( x )  &=     -  \lambda c  ( x )  \nabla_y w ( x, y ),  &  S \times \Omega \times Y_p  ,
  \label{2.3a}
  \\
  \nabla_y \cdot \vec{u} (x,y) &= 0 , &  S \times \Omega \times Y_p ,
  \label{2.3b}
  \\
  \nabla_x \cdot \overline{ \vec{u} } ( x )
  & = 0 ,  &   S \times \Omega ,
  \label{2.3c}
  \\
  \vec{u} (x,y) & = 0  ,  &     S \times \Omega \times \Gamma_s  ,
  \label{2.3d}
  \\
  \partial_t c(x) &= \Delta_y w (x,y)  ,  &    S \times \Omega \times Y_p ,
  \label{2.3e}
  \\
  \overline{w} ( x )   + \nabla_x \cdot  \overline{\nabla_y c_1} ( x ) &=  -  \Delta_x c (x) + f ( c(x) )   , &   S \times \Omega ,  
  \label{2.3f}
  \\
  \nabla_y \cdot \nabla_y c_1 (x,y) &= 0  , &  S \times \Omega \times Y_p ,
  \label{2.3g} 
  \\
  c ( 0, x ) &= c_0 (x)  , &  \Omega.
  \label{2.3h}
 \end{align}
\end{subequations}
where $ \bar{\kappa} ( x ) = \frac{1}{|Y_p|}  \int_{ \partial Y_p} \kappa (x,y) \,dy $, $ x \in \Omega$ denotes the mean of the quantity $\kappa$ over the pore space $Y_p$. Also, $c_1(t,x,y) = \psi(t,x) \xi(y)$, where $\xi (y)$ is a linear function. From \eqref{2.3e}, $ w( x,y ) = \partial_t c(x) ~ \varsigma ( y ) $, and hence the cell problem is as follows: 
\begin{equation}
\begin{cases}
&  \partial_{y_i y_j} \varsigma ( y ) = \delta_{ij} \quad \text{ in } Y_p, 
\\
&  \vec{n} \cdot \nabla_y \varsigma ( y ) = 0  \quad \text{ on }   \partial Y_p , 
\\
&  \varsigma ( y ) \text{ is $Y_p$-periodic}. 
\end{cases}
\end{equation}

The systems of equations  \eqref{2.3a}-\eqref{2.3h} is the required homogenized (upscaled) model of \eqref{1.5a}-\eqref{1.5h}.

\end{theorem}


\section{Anticipated upscaled model via Asymptotic expansion method}
In this section we derive the homogenized version of the model $\mathcal{(P^\eps)}$ as described in the previous section via some formal method, namely asymptotic expansion, which does not talk about any convergence. Now as per asymptotic expansion technique let us consider the following expansions:	 
\begin{align} \label{3.1}
		\vec {u^{\varepsilon}}  &= \sum_{i=0}^\infty\eps^i \vec {u_i}, c^{\varepsilon}   = \sum_{i=0}^\infty \eps^ic_i, w^{\varepsilon}  = \sum_{i=0}^\infty \eps^iw_i \text{ and } p^{\varepsilon}   = \sum_{i=0}^\infty \eps^ip_i, 
\end{align}
where each term $ \vec {u_i}$ or $p_i$ or $c_i$ and $w_i$ are $Y$-periodic functions in $y$. We know that $\nabla=\nabla_x+\frac{1}{\varepsilon}\nabla_y$. We substitute the expressions for $\vec{u}^\eps, c^\eps, w^\eps, p^\eps$ in the problem $\mathcal{(P^\eps)}$, then 
\begin{align*}
	-  \mu \varepsilon^2 \big(  \varepsilon ^{-2} \Delta_y + \varepsilon ^{-1} (\nabla_x \cdot \nabla_y & + \nabla_y \cdot \nabla_x) +  \varepsilon ^0 \Delta_x \big) ~ \big(\vec {u_0}   + \varepsilon \vec {u_1}  + \varepsilon ^2 \vec {u_2} +...\big) 
	\\
	&+ \big( \nabla_x +  \varepsilon^{-1} \nabla_y \big) ~ \big( p_0   
	+ \varepsilon p_1  + \varepsilon ^2 p_2...\big) 
	\\
	&= -\lambda \varepsilon ~ \big( c_0  + \varepsilon c_1   + \varepsilon ^2 c_2... \big) ~ \big( \nabla_x +  \varepsilon^{-1} \nabla_y \big) ~ \big(w_0  + \varepsilon w_1   + \varepsilon ^2 w_2... \big)
	 \\
	 \implies - \mu ( \varepsilon ^0 \Delta_y \vec{u_0} + \varepsilon   ( \Delta_y \vec{u_1} &+ (\nabla_x \cdot \nabla_y + \nabla_y \cdot \nabla_x) \vec{u_0} ) + \varepsilon ^2 (...) ) 
	 \\
	 &+ \varepsilon ^{-1} \nabla_y p_0 + \varepsilon ^0 (\nabla_x p_0 
		   \nabla_y p_1) 
		+ \varepsilon ( \nabla_x p_1 + \nabla_y p_2 ) + \varepsilon ^2 (...)
	\\
	&= - \lambda ( \varepsilon ^0 (c_0 \nabla_y w_0) + \varepsilon   ( c_1 \nabla_y w_0 + c_0 ( \nabla_x w_0 + \nabla_y w_1 ) ) + \varepsilon ^2 (...) ) 
	\end{align*}
	 \begin{gather} \label{3.2}
	 \begin{aligned}
	 \implies   \varepsilon ^{-1} (\nabla_y p_0)  + \varepsilon ^0 ( - \mu \Delta_y \vec{u_0} &+  \nabla_x p_0 + \nabla_y p_1 ) + \varepsilon  ( - \mu( \Delta_y \vec{u_1} 
	 \\
	 &+ (\nabla_x \cdot \nabla_y + \nabla_y \cdot \nabla_x) \vec{u_0} ) 
	 +  \nabla_x p_1 
	 + \nabla_y p_2  )  + \mathcal{O} (\varepsilon) 
	 \\
	 &= \varepsilon ^0 ( -\lambda (c_0 \nabla_y w_0) ) + \varepsilon   (-\lambda ( c_1 \nabla_y w_0 + c_0 ( \nabla_x w_0 + \nabla_y w_1 ) ) ).
	 \end{aligned}
	 \end{gather}
	 
We use \eqref{3.1} in \eqref{1.5b} then
\begin{equation} \label{3.3}
		\varepsilon ^{-1} \nabla_y \cdot \vec{u_0} + \varepsilon ^0 ( \nabla_x \cdot \vec{u_0} + \nabla_y \cdot \vec{u_1} ) + \varepsilon ( \nabla_x \cdot \vec{u_1} + \nabla_y \cdot \vec{u_2} ) + \varepsilon ^2 (...) = 0.
\end{equation}
	
From \eqref{1.5d}, after plugging the expansions, we obtain 
\begin{align*}
		\partial_t (  c_0   + \varepsilon c_1   + \varepsilon ^2 (...) ) &+ \varepsilon ( \varepsilon ^{-1} \nabla_y + \nabla_x )   \cdot \{ ( c_0   + \varepsilon c_1   + \varepsilon ^2 (...) )  (   \vec {u_0}   + \varepsilon \vec {u_1}   + \varepsilon ^2(...) ) \} 
		\\
		&= \varepsilon^2 ( \varepsilon ^{-2} \Delta_y w_0 
		+ \varepsilon ^{-1} (  \Delta_y w_1 + ( \nabla_x \cdot \nabla_y + \nabla_y \cdot \nabla_x ) w_0  ) 
		\\
		&+ \varepsilon ^0 ( \Delta_y w_2 + ( \nabla_x \cdot \nabla_y + \nabla_y \cdot \nabla_x  ) w_1 + \Delta_x w_0 ) )  
\end{align*}  
\begin{align}
	 	\implies \partial_t (  c_0   + \varepsilon c_1 ) +  \varepsilon ^0 \{\nabla_y   \cdot ( c_0 \vec{u_0} ) \} &
	 +\varepsilon \{  \nabla_y   \cdot ( c_0 \vec {u_1} ) +  \nabla_x   \cdot ( c_0  \vec{u_0}  ) + \nabla_y \cdot (  c_1 \vec{u_0} ) \}   
	 \notag\\
	 &=   \varepsilon ^0 \Delta_y w_0 
	 + \varepsilon ^1 (  \Delta_y w_1 + ( \nabla_x \cdot \nabla_y + \nabla_y \cdot \nabla_x ) w_0  ) 
	 \notag\\
	 &+ \varepsilon ^2 ( \Delta_y w_2 + ( \nabla_x \cdot \nabla_y + \nabla_y \cdot \nabla_x  ) w_1 + \Delta_x w_0 ) + \mathcal{O} (\varepsilon).
	 \label{3.4}
\end{align} 

Next, we substitute the expansions for $w_{\varepsilon}$, $c_{\varepsilon}$ in \eqref{1.5e} and use the Taylor series expansion of $f$ around $c_0$ which leads to
\begin{align}
	w_0  + \varepsilon w_1   &= \varepsilon ^{-2}( - \Delta_y c_0 ) + \varepsilon ^{-1} ( - \Delta_y c_1 - ( \nabla_x \cdot \nabla_y + \nabla_y \cdot \nabla_x ) c_0 ) 
	\notag\\
	& + \varepsilon ^0 ( - \Delta_y c_2 - ( \nabla_x \cdot \nabla_y + \nabla_y \cdot \nabla_x ) c_1 - \Delta_x c_0  + f( c_0 ) ) 
	\notag\\
	& + \varepsilon (   - \Delta_y c_3 - ( \nabla_x \cdot \nabla_y + \nabla_y \cdot \nabla_x ) c_2 - \Delta_x c_1   +  c_1 \frac{f(c_0)}{c_0} )  +  \mathcal{O} (\varepsilon).
	\label{3.5} 
\end{align} 

Now we substitute the expansions in the boundary conditions. From \eqref{1.5c}, we obtain
\begin{equation} \label{3.6}
  \vec {u_0}   + \varepsilon \vec {u_1}  + \varepsilon ^2 \vec {u_2} +... = 0 \quad \mbox{on }(0,T) \times \partial\Omega_p^ \varepsilon.
\end{equation}

From \eqref{1.5f} and \eqref{1.5g}, we get 
\begin{equation*}
( \varepsilon ^{-1} \nabla_y + \nabla_x )( c_0   + \varepsilon c_1   + \varepsilon ^2 c_2 + ... ) \cdot \vec{n} = 0 
\end{equation*}
\begin{equation} \label{3.7}
\implies \varepsilon ^{-1} \nabla_y c_0 \cdot \vec{n} + \varepsilon^0 ( \nabla_x c_0 + \nabla_y c_1 ) \cdot \vec{n} + \varepsilon ( \nabla_x c_1 + \nabla_y c_2 ) \cdot \vec{n} + ... = 0
\end{equation}
and
\begin{equation*}
( \varepsilon ^{-1} \nabla_y + \nabla_x )( w_0   + \varepsilon w_1   + \varepsilon ^2 w_2 + ... ) \cdot \vec{n} = 0 
\end{equation*}
\begin{equation} \label{3.8}
\implies  \varepsilon ^{-1} \nabla_y w_0 \cdot \vec{n} + \varepsilon^0 ( \nabla_x w_0 + \nabla_y w_1 ) \cdot \vec{n} + \varepsilon ( \nabla_x w_1 + \nabla_y w_2 ) \cdot \vec{n} + ... = 0
\end{equation}
respectively.

We equate the coefficients of $\varepsilon ^{-1}$ from \eqref{3.2} and coefficient of $\varepsilon ^ {-2}$ from \eqref{3.5}, respectively, then	
\begin{align} \label{3.9}
		\nabla_y p_0 &= 0,   &\text{and}     &        &  \Delta_y c_0 = 0 \quad \text{for } y \in Y_p.
\end{align} 
and from \eqref{3.9} it follows that
\begin{align} \label{3.14}
		p_0 &= p_0(x),  \quad \text{ for } y \in Y_p.
\end{align}
	
Since $ \nabla_x p_0 = \sum_{i=1}^{n} e_i \pdv{p_0}{x_i} $, the coefficient of $\varepsilon ^0$  from \eqref{3.2} gives 
	\begin{align} 
	&- \mu \Delta_y \vec{u_0} + \sum_{i=1}^{n} e_i \pdv{p_0}{x_i}  + \nabla_y p_1 = -\lambda  c_0 \nabla_y w_0
		\notag\\
	\implies &     - \mu \Delta_y \vec{u_0} + \nabla_x p_0  + \nabla_y p_1 = \lambda c_0 \nabla_y w_0 . 
		\label{3.15}
	\end{align}
 
Again, upon equating $\varepsilon ^{-1}$ and $ \varepsilon ^ 0 $ coefficients of \eqref{3.3}, we get respectively 
\begin{align} 
 \nabla_y \cdot \vec{u_0} = 0  \quad \text{in } S \times \Omega \times Y_p,
 \label{3.18a}
 \\
  \nabla_x \cdot \vec{u_0} + \nabla_y \cdot \vec{u_1}  = 0  \quad \text{in } S \times \Omega \times Y_p.
   \label{3.18b}
\end{align}

Equating $\varepsilon ^0$, $\varepsilon $ coefficients of \eqref{3.6} we get 
\begin{align} 
\vec{u_0} &= 0,   & &  \vec{u_1} = 0 \quad \text{on } S \times \Omega \times \Gamma_s.
\label{3.19}
\end{align}

We integrate \eqref{3.18b} over $Y_p$, which leads to
\begin{align}  
&\int_{Y_p} \{ \nabla_x \cdot \vec{u_0} + \nabla_y \cdot \vec{u_1} \} \,dy = 0
\notag \\
\implies  & \nabla_x    \cdot  \int_{Y_p}  \vec{u_0}  \,dy  +  \int_{ \partial Y_p}  \vec{n} \cdot \vec{u_1} \,d \sigma (y)   = 0  .
\label{3.25}
\end{align}

Using \eqref{3.19} the boundary integral in \eqref{3.25} vanishes and we obtain
\begin{equation}
 \nabla_x    \cdot  \int_{Y_p}  \vec{u_0} ( x,y )  \,dy  =  0   \quad \text{  in $ S \times \Omega  $.} 
\end{equation}

The coefficient of $\varepsilon ^{-1}$ from \eqref{3.5} gives using \eqref{3.14}
\begin{align} 
& \Delta_y c_1 = - ( \nabla_x \cdot \nabla_y + \nabla_y \cdot \nabla_x ) c_0 , 
\label{3.14a} \\
\implies & \int_{\partial Y_p} \vec{n} \cdot  ( \nabla_y c_1  + \nabla_x   c_0 )  \,d \sigma (y)  = -  \nabla_x \cdot \Big\{\frac{1}{|Y_p|} \int_{Y_p} \nabla_y c_0 \,dy \Big\},
\label{3.24}
\end{align} 

From \eqref{3.7}, the coefficient of $\varepsilon ^0$ gives
\begin{equation} \label{3.34} 
\vec{n} \cdot  ( \nabla_y c_1  + \nabla_x   c_0 ) = 0. 
\end{equation}

Using \eqref{3.34} in \eqref{3.24} leads to
\begin{equation} \label{3.35} 
\nabla_x \cdot \Big\{\frac{1}{|Y_p|} \int_{Y_p} \nabla_y c_0 \,dy \Big\} = 0.
\end{equation} 

From \eqref{3.35} and \eqref{3.9}, clearly $ c_0  = c_0 ( t, x )$.	Also from \eqref{1.5h} we get, $ c_{0} ( 0,x ) = c_0 ( x ) $. 

From \eqref{3.4}, the coefficient of $\varepsilon ^0$ gives
	\begin{align}  
	  	\partial_t c_0 + \nabla_y \cdot ( c_0 \vec{u_0} )  &= \Delta_y w_0  
		 \label{3.20}
	\end{align}
	
Using \eqref{3.18a} we get from \eqref{3.20}
\begin{align}  
\partial_t c_0  &= \Delta_y w_0 . 
\label{3.21}
\end{align}

By the separation of variables we obtain
\begin{align} \label{3.22}
 w_0 ( t,x,y ) =  \partial_t c_0 (t,x ) \varsigma ( y )  .
\end{align}

We equate coefficients of $ \varepsilon ^{-1} $ from \eqref{3.8} and obtain after combining it with (3.22).
\begin{subequations} \label{3.23}
\begin{align}
\partial_{y_i y_j} \varsigma ( y ) = \delta_{ij}, \quad \text{ in } Y_p
\label{3.23a}
 \\
\vec{n} \cdot\nabla_y \varsigma (y) = 0, \quad \text{on } \partial Y_p
\label{3.23b}
 \\
\varsigma (y) \text{ is $Y_p$-periodic. }
\label{3.23c}
\end{align}
\end{subequations}

From \eqref{3.14a} we get 
\begin{align}
 \Delta_y c_1 ( t,x,y )  &= 0, \quad \text{in } S \times \Omega \times Y_p .
 \label{3.26a}
\end{align}

Thus $ c_1 $ has the form $ c_1 ( t,x,y ) = \psi ( t,x ) \xi ( y )  $, where $\xi ( y )$ is a linear function in variable $y$.

The coefficient of $\varepsilon ^0$ from \eqref{3.5} yields 
	\begin{align}
		w_0 &= - ( \Delta_y c_2 + (  \nabla_x \cdot \nabla_y + \nabla_y \cdot \nabla_x  ) c_1 + \Delta_x c_0 ) + f( c_0 ), \notag\\
		\implies \Delta_y c_2  &=  -  w_0  - (  \nabla_x \cdot \nabla_y + \nabla_y \cdot \nabla_x  ) c_1 - \Delta_x c_0   + f( c_0 ), \label{3.28}
	\end{align}

Next we integrate \eqref{3.28} over $Y_p$.
\begin{align} 
 \implies & \int_{Y_p} \Delta_y c_2 \,dy =   \int_{Y_p} \{-  w_0  - (  \nabla_x \cdot \nabla_y + \nabla_y \cdot \nabla_x  ) c_1 - \Delta_x c_0   + f( c_0 ) \} \,dy , 
  \notag
  \\   
 \implies & \int_{\partial Y_p} \vec{n} \cdot  ( \nabla_y c_2  + \nabla_x   c_1 )  \,d \sigma (y) =  \int_{Y_p} \{-  w_0  -  \nabla_x \cdot \nabla_y  c_1 - \Delta_x c_0   + f( c_0 ) \} \,dy
 \label{3.31} 
\end{align}

From \eqref{3.7}, the coefficient of $\varepsilon$ gives
\begin{equation} \label{3.32} 
  \vec{n} \cdot  ( \nabla_y c_2  + \nabla_x   c_1 ) = 0 
\end{equation}

The left hand side of \eqref{3.31} vanishes using \eqref{3.32}. Thus we get from \eqref{3.31}
\begin{align} 
  w_0 ( t,x,y )   +  \Delta_x c_0 ( t,x  )  +  \nabla_x \cdot \nabla_y  c_1  ( t,x,y )   =   f( c_0 ( t,x  ) ) \quad \text{in } S \times \Omega \times  Y_p  
\label{3.33} 
\end{align} 

From \eqref{3.15} and \eqref{3.22}, we have
\begin{align} 
 - \mu \Delta_y \vec{u_0}   + \nabla_y p_1 =  - \nabla_x p_0  -  \lambda  ~ c_0 ~  \partial_t c_0 ~  \nabla_y \varsigma ( y )
\notag \\
- \mu \Delta_y \vec{u_0}   + \nabla_y p_1 =  - \nabla_x p_0 ( x ) -  \frac{\lambda}{2}  \frac{d}{dt}  \{c_0 (t,x) \}^2 ~  \nabla_y \varsigma ( y ) 
\label{3.36}
\end{align}


\section{Proof of theorem 2.1}

\subsection{A priori Estimates}
We now proceed to derive some a priori estimates to show that the sequences of functions $\vec{u^{\varepsilon}}, c^ \varepsilon $ and $ w^ \varepsilon$
are bounded independently of $\varepsilon$ in appropriate function spaces. We set $\eta = \vec{u^{\varepsilon}} $, $\phi =   \lambda w^{\varepsilon}$ and $\psi =  \lambda \partial_t c^{\varepsilon}$ and use $\nabla (c^{\varepsilon} w^{\varepsilon}) = c^{\varepsilon} \nabla w^{\varepsilon} + w^{\varepsilon} \nabla c^{\varepsilon}$. From \eqref{2.1} we obtain
\begin{align*}  
	&	\mu \varepsilon ^2   \int_{S \times \Omega_p^ \varepsilon} | \nabla \vec{u^\varepsilon} |^2 \,dx \,dt +   \varepsilon  \lambda   \int_{S \times \Omega_p^ \varepsilon} c^{\varepsilon} \nabla w^{\varepsilon} \cdot \vec{u^{\varepsilon}} \,dx \,dt = 0,
	\\
	&	\lambda \int_{S} \langle \partial_t c^{\varepsilon} , w^{\varepsilon} \rangle \,dt -   \lambda \varepsilon  \int_{S \times \Omega_p^ \varepsilon} c^{\varepsilon} \vec{u^{\varepsilon}} \cdot \nabla w^{\varepsilon} \,dx \,dt +    \lambda \varepsilon^2  \int_{S \times \Omega_p^ \varepsilon} | \nabla
		w^{\varepsilon} | ^2  \,dx \,dt = 0, 
	\\
		-  &  \lambda \int_{S} \langle \partial_t c^{\varepsilon} ,  w^{\varepsilon} \rangle  \,dt +    \frac{  \lambda}{2} \frac{d}{dt} \int_{\Omega_p^ \varepsilon} | \nabla c^{\varepsilon} |^2 \,dx +    \lambda \int_{S} \langle  f(c^{\varepsilon}) , \partial_t c^{\varepsilon} \rangle   \,dt = 0.
\end{align*}
This implies
\begin{align} 
	&	\mu \varepsilon ^ {2 }  \int_{S \times \Omega_p^ \varepsilon}  | \nabla \vec{u^\varepsilon} |^2 \,dx \,dt  + \lambda \varepsilon  \int_{S \times \Omega_p^ \varepsilon} c^{\varepsilon} \nabla w^{\varepsilon} \cdot \vec{u^{\varepsilon}} \,dx \,dt = 0,
	\label{4.1}\\
	&	\lambda \int_{S} \langle \partial_t c^{\varepsilon} , w^{\varepsilon} \rangle \,dt  - \lambda \varepsilon   \int_{ S \times \Omega_p^ \varepsilon} c^{\varepsilon} \vec{u^{\varepsilon}} \cdot \nabla w^{\varepsilon} \,dx \,dt + \lambda \varepsilon^2   \int_{ S \times \Omega_p^ \varepsilon} | \nabla
		w^{\varepsilon} | ^2  \,dx \,dt = 0, 
	\label{4.2}\\
	- &   \lambda  \int_{S} \langle \partial_t c^{\varepsilon} , w^{\varepsilon} \rangle \,dt  +  \frac{  \lambda}{2} \frac{d}{dt} \int_{\Omega_p^ \varepsilon} | \nabla c^{\varepsilon} |^2 \,dx  +   \lambda \frac{d}{dt} \int_{\Omega_p^ \varepsilon} F(c^{\varepsilon}) \,dx = 0, 
	\label{4.3}
\end{align}
where $F$ is considered to be a quadratic double-well free energy functional as
\begin{equation} \label{4.4}
	F(c^{\varepsilon}) = \frac{1}{4} (({c^{\varepsilon}})^2 - 1)^2.
\end{equation}
We add \eqref{4.1}-\eqref{4.3} and integrate over $(0,t)$, then
\begin{align} 
		\mu \varepsilon ^ {2 } \int_{0}^{t} \int_{\Omega_p^ \varepsilon}  | \nabla \vec{u^\varepsilon} |^2 \,dx \,dt   &+  \lambda \varepsilon ^2 \int_{0}^{t} \int_{  \Omega_p^ \varepsilon} | \nabla w^{\varepsilon} | ^2  \,dx \,dt   +  \frac{  \lambda}{2} \int_{\Omega_p^ \varepsilon} | \nabla c^{\varepsilon}(t) |^2 \,dx  
		+    \lambda \int_{\Omega_p^ \varepsilon} F(c^{\varepsilon}(t)) \,dx 
		\notag\\
		&= \frac{  \lambda}{2} \int_{\Omega_p^ \varepsilon} | \nabla c^{\varepsilon}_0 (x) |^2 \,dx  +   \lambda \int_{\Omega_p^ \varepsilon} F(c^{\varepsilon}_0 (x)) \,dx. 
		\label{4.5}
\end{align}
We note that $c_0\in L^\infty(\Omega)\cap H^1(\Omega)$ such that $ \sup_{\eps>0}||c_0||_{L^\infty(\Omega)\cap H^1(\Omega)}<\infty$ and $F(c^{\varepsilon}) \geq 0$. Therefore, \eqref{4.5} implies 
\begin{align}   \label{4.6}
		 \sqrt{\mu} \varepsilon   \norm{\nabla \vec{u^\varepsilon} }_{L^2(S \times {\Omega}_p^{\varepsilon})}   +  \sqrt{ \lambda} \varepsilon \norm{ \nabla w^{\varepsilon} }_{L^2(S \times {\Omega}_p^{\varepsilon})} + \sqrt{\frac{\lambda}{2}} \norm{ \nabla c^{\varepsilon} }_{L^{\infty} (S ; L^2 ({\Omega}_p^{\varepsilon}) )}  \leq C  
\end{align}
Next, by Young's inequality, from \eqref{4.5} we obtain
\begin{equation} \label{4.10}
		\int_{\Omega_p^ \varepsilon} F(c^{\varepsilon}(t)) \,dx =\frac{1}{4}  \int_{  \Omega_p^ \varepsilon} ((c^{\varepsilon})^2 - 1)^2 \,dx \leq C\Rightarrow\int_{  \Omega_p^ \varepsilon} | c^{\varepsilon} |^4 \,dx \leq C \; \forall t\Rightarrow \sup_{\varepsilon>0}\norm{ c^{\varepsilon} }_{L^{\infty}(S ; L^4 ({\Omega}_p^{\varepsilon}))} \leq C.
\end{equation} 
We set $\psi = 1$ as a test function in \eqref{1.5e}, then
\begin{equation} \label{4.8}
		\int_{  \Omega_p^ \varepsilon} w^{\varepsilon} \,dx = \int_{  \Omega_p^ \varepsilon} f(c^{\varepsilon}) \,dx\Rightarrow \left|\int_{  \Omega_p^ \varepsilon} w^{\varepsilon} \,dx\right|\le \int_{  \Omega_p^ \varepsilon}| f(c^{\varepsilon}) |\,dx\le\int_{  \Omega_p^ \varepsilon}(|(c^{\varepsilon}) |^3+|(c^{\varepsilon}) |)\,dx \le C \text{ by \eqref{4.10}}.
\end{equation}	
By Poincare's inequality, we have
\begin{equation} \label{4.9}
	||w^\varepsilon -\int_{\Omega^\varepsilon_p}w^\varepsilon\,dx||_{L^2(\Omega^\varepsilon_p)}\le C||\nabla w^\varepsilon||_{L^2(\Omega^\varepsilon_p)}	\Rightarrow \norm{ w^{\varepsilon} }_{L^2(S \times {\Omega}_p^{\varepsilon})} \leq C.
\end{equation}  

By Gagliardo-Nirenberg-Sobolev inequality for Lipschitz domain, we obtain $||u^\varepsilon||_{L^4(Y)}\le C||\nabla u^\varepsilon||_{L^2(Y)}$, where $C$ depend on $n$ and $Y$. By imbedding theorem, we have $||u^\varepsilon||_{L^2(Y)}\le C||u^\varepsilon||_{L^4(Y)}\le C$. By a straightforward scaling argument, we obtain
\begin{equation}
	\norm{\vec {u^{\varepsilon}}}_{L^4(\Omega_p^ \varepsilon)} \leq C.	\label{4.17}
\end{equation}
From \eqref{4.2} we get, 
\begin{align}
		&\int_{S} \langle \partial_t c^{\varepsilon} , \phi \rangle \,dt  = \varepsilon  \int_{S \times \Omega_p^ \varepsilon} c^{\varepsilon} \vec{u^{\varepsilon}} \cdot \nabla \phi \,dx \,dt - \varepsilon ^2  \int_{S \times \Omega_p^ \varepsilon} \nabla
		w^{\varepsilon} \cdot \nabla \phi \,dx \,dt ,
		\notag\\
		\implies& | \langle \partial_t c^{\varepsilon} , \phi \rangle | \leq \varepsilon   \norm{c^{\varepsilon}}_{L^4(\Omega_p^ \varepsilon)} \norm{\vec {u^{\varepsilon}}}_{L^4(\Omega_p^ \varepsilon)} \norm{\nabla \phi}_{L^2(\Omega_p^ \varepsilon)} + \varepsilon ^2 \norm{\nabla w^{\varepsilon}}_{L^2(\Omega_p^ \varepsilon)}\norm{\nabla \phi}_{L^2(\Omega_p^ \varepsilon)}
        \notag\\
		\implies& \sup _{\norm{\phi}_{\mathbf{H}^1 ({\Omega}_p ^ \varepsilon)} \leq 1} | \langle \partial_t c^{\varepsilon} , \phi \rangle |  \leq    \norm{c^{\varepsilon}}_{L^4(\Omega_p^ \varepsilon)} \norm{\vec {u^{\varepsilon}}}_{L^4(\Omega_p^ \varepsilon)}  + \varepsilon  \norm{\nabla w^{\varepsilon}}_{L^2(\Omega_p^ \varepsilon)}\notag\\
		\implies&  \norm{\partial_t c^{\varepsilon}}_{L^2( S ; H^1 ({\Omega}_p ^ \varepsilon)^* )} \leq C ~~~ \quad  \forall \varepsilon > 0\label{4.16}
\end{align}
	
	
From proposition III.1.1 in \cite{temam2001navier} and \eqref{2.1a}, there exist a pressure $  p^{\varepsilon} \coloneqq \partial_t P^{\varepsilon} \in W^{-1, \infty} ( S ,$ $ L_0^2 ( \Omega_p^{\varepsilon} ) )  $ such that
\begin{align*}
& \mu \varepsilon ^2 \int_{S \times \Omega_p^ \varepsilon} \nabla \vec{u^\varepsilon} : \nabla \eta \,dx \,dt -\int_{\Omega_p^ \varepsilon} P^{\varepsilon} (t) \nabla \cdot \eta \,dx=-\varepsilon  \lambda \int_{S \times \Omega_p^ \varepsilon} c^{\varepsilon} \nabla w^{\varepsilon} \cdot \eta \,dx \,dt\quad  \forall \eta \in H_0^1( \Omega_p^{\varepsilon} )^n,\notag\\
\Rightarrow& \langle \nabla P^{\varepsilon} (t) , \eta \rangle \leq \mu \varepsilon^2  \int_S \norm{ \nabla \vec{u}^{\varepsilon}}_{L^2(\Omega_p^{\varepsilon})} \norm{ \nabla \eta^{\varepsilon}}_{L^2(\Omega_p^{\varepsilon})} \,dt + \int_S \norm{ c^{\varepsilon} }_{L^4(\Omega_p^{\varepsilon})} \varepsilon \norm{ \nabla w^{\varepsilon}}_{L^2(\Omega_p^{\varepsilon})} \,dt.
\end{align*}
Thus by \eqref{4.6} and \eqref{4.10} it immediately follows that
\begin{align}
 \langle \nabla P^{\varepsilon} (t) , \eta \rangle 
 \leq
  C \norm{ \eta }_{H_0^1( \Omega_p^{\varepsilon} )^n}  \Rightarrow
\sup_{t \in [0,T]}  \norm { \nabla P^{\varepsilon} (t) }_{H^{-1}( \Omega_p^{\varepsilon} )^n}  
 \leq
  C       \quad \forall \varepsilon > 0   .
\end{align}
Now, with the help of a-priori estimates from (2.3), the existence of solution of $\mathcal{(P^\eps)}$ can be shown using the similar arguments from \cite{Fen06,SL03,FHC07}.

\section{Proof of Theorem 2.2 (Homogenization of problem $\mathcal{(P^\eps)}$)}
 
In the previous section, we established the existence of a weak solution of the problem $\mathcal{(P^\eps)}$. In this section we rigorously derive the upscaled model for $\varepsilon \rightarrow 0$. Two-scale convergence introduced in section 2, is a type of weak convergence in some $L^p$-space. The idea of two-scale convergence resides on the assumption that the oscillating sequence of functions are defined over some fixed domain, say $\Omega $ and we have boundedness of such sequence of functions in $L^p(\Omega)$ for some $p$. Since there is no oscillation in $t\in S$, we are only focused on the oscillation in $x\in \Omega$. As our solutions $c^{\varepsilon}$, $w^{\varepsilon}$, $\vec{u^{\varepsilon}}$ and $P^{\varepsilon}$ are defined over $\Omega_p^\varepsilon$, in order to apply two-scale convergence we need to obtain the a-priori estimates for $c^{\varepsilon}$, $w^{\varepsilon}$, $\vec{u^{\varepsilon}}$, $P^{\varepsilon}$ etc. in $L^p(S\times\Omega)$ for some $p$. This is not so straightforward. Usually, we first obtain the estimates for $c^{\varepsilon}$, $w^{\varepsilon}$, $\vec{u^{\varepsilon}}$, $P^{\varepsilon}$ in $\Omega_p^\varepsilon $ and then use the extension operator defined in lemma 2.2 to extend these estimates to all of $S\times\Omega$. We start with the construction of an extension of solution from $\Omega_p^{\varepsilon}$ to $\Omega$ in the lemma below. 


\begin{lemma} 
	  There exists a positive constant $C$ depending on $c_0$, $\vec{u_0}$, $n$, $|Y|$, $\lambda$ and $\mu$ but independent of $\varepsilon$ and extensions ($\tilde{c^{\varepsilon}}$, $\tilde{w^{\varepsilon}}$, $\tilde{\vec{u}}^{\varepsilon}$, $\tilde{P^{\varepsilon}}$) of the solution ($c^{\varepsilon}$, $w^{\varepsilon}$, $\vec{u^{\varepsilon}}$, $P^{\varepsilon}$) to $S \times \Omega$ such that 
	   \begin{small}
	   	\begin{align} 
	   	&\norm{\tilde{\vec{u}}^{\varepsilon}}_{L^{\infty}( S ;L^2(\Omega)^n)} +
	   	\sqrt{\mu} \varepsilon \norm{\nabla \tilde{\vec{u}}^{\varepsilon}}_{L^2( S \times \Omega )^{n \times n}} + \norm{\tilde{w}^{\varepsilon}}_{L^2(S;H^1(\Omega))} +  
	   	\sqrt{\lambda } \varepsilon \norm{ \nabla \tilde{w}^{\varepsilon} }_{L^2(S \times \Omega)^n}
	   	\notag\\
	   	&+\norm{\tilde{c}^{\varepsilon}}_{L^{\infty}(S;L^4(\Omega))} + \sqrt{\frac{\lambda}{2}} \norm{\nabla \tilde{c}^{\varepsilon}}_{L^{\infty}(S;L^2(\Omega)^n)} + 
	   	\norm{ \partial_t \tilde{c}^{\varepsilon}}_{L^2(S;H^1(\Omega)^{*})} +
	   	\sup_{t \in [0,T]} \norm{\tilde{P}^{\varepsilon}(t)}_{L_0^2(\Omega)}  
	   	\leq C.
	   	\label{5.1}	
	   	\end{align}
	   \end{small}
	    
\end{lemma}

\begin{proof}
The extensions of $\tilde{c^{\varepsilon}}$, $\tilde{w^{\varepsilon}}$, $\tilde{\vec{u}}^{\varepsilon}$ to all of $ S \times \Omega$ follows from the lemma 2.2 and estimate (2.2). Next, we consider the extension of $\partial_t c^{\varepsilon}$ from $L^2(S;H^1(\Omega^\varepsilon_p)^*)$ to $L^2(S;H^1(\Omega)^*)$. For $\Theta \in H^1(\Omega^\varepsilon_p)^*$, we define the extension operator $F^{\varepsilon} : H^1(\Omega ^{\varepsilon}_p)^* \rightarrow H^1(\Omega)^*$ as
	 \begin{equation} \label{5.2}
	  \langle F^{\varepsilon} \Theta , \phi \rangle_{H^1(\Omega)^* \times H^1(\Omega)} = \langle \Theta , \mathcal{R}^{\varepsilon} \phi \rangle_{H^1(\Omega^{\varepsilon}_p)^* \times H^1(\Omega^{\varepsilon}_p)},
	 \end{equation}  
	 where $\mathcal{R}^{\varepsilon} : H^1(\Omega) \rightarrow H^1(\Omega^{\varepsilon}_p)$ is the restriction operator $\mathcal{R}^{\varepsilon} \phi = \phi |_{\Omega^{\varepsilon}_p}$ for $\phi \in H^1(\Omega)$. Since $\norm{  {R}^{\varepsilon} \phi }_{H^1(\Omega)} \leq C\norm{\phi}_{H^1(\Omega^{\varepsilon}_p)}$ it follows that 
	 \begin{equation} \label{5.3}
	 \norm{F^{\varepsilon} \Theta}_{H^1(\Omega)^*} \leq C\norm{\Theta}_{H^1(\Omega^{\varepsilon}_p)^*}.
	 \end{equation}
Using \eqref{5.2} we define the extension $\widetilde{\partial_t c^{\varepsilon}}$ of $\partial_t c^{\varepsilon}$ in $L^2(S;H^1(\Omega)^*)$ as
	 \begin{align} \label{5.4}
	  \int_{0}^{T} \langle \widetilde{\partial_t c^{\varepsilon}} , \phi \rangle_{H^1(\Omega)^* \times H^1(\Omega)} \coloneqq \int_{0}^{T} \langle F^{\varepsilon} \partial_t c^{\varepsilon} , \phi \rangle_{H^1(\Omega)^* \times H^1(\Omega)},\notag\\
	  \Rightarrow ||\widetilde{\partial_t c^{\varepsilon}}||_{L^2(S;H^1(\Omega)^*)}=||F^{\varepsilon} \partial_t c^{\varepsilon}||_{L^2(S;H^1(\Omega)^*)}\le C|| \partial_t c^{\varepsilon}||_{L^2(S;H^1(\Omega_p^\varepsilon)^*)}\le C<\infty.
	 \end{align}
Following the arguments as in \cite{All93,bavnas2017homogenization}, the pressure $P^{\varepsilon}$ can be extended into whole $\Omega$, i.e.
	 \begin{equation} \label{5.5}
	  \sup_{t \in [0,T]} \norm{ \tilde{P}^{\varepsilon}(t) }_{L_0^2(\Omega)} \leq C,
	 \end{equation}
	 where $C$ is independent of $\varepsilon$.
\end{proof} 


\begin{lemma}
	 Let ($\vec{u^{\varepsilon}}$, $P^{\varepsilon}$, $c^{\varepsilon}$, $w^{\varepsilon}$)$_{\varepsilon > 0}$ be the extension of the weak solution from Lemma 5.1 (denoted by the same symbol). Then there exists some functions $\vec{u} \in L^2( S \times \Omega ; H^1_{\#}(Y) )^n$, $c \in L^2( S ;H^1(\Omega) )$, $P \in L^2( S \times \Omega \times Y )$, $c_1$, $w \in L^2( S \times \Omega; H^1_{\#}(Y) )$ and a subsequence of $ ( \vec{u^{\varepsilon}}, P^{\varepsilon}, c^{\varepsilon}, w^{\varepsilon} )_{\varepsilon > 0} $, still denoted by the same symbol, such that the following convergences hold:
	 \begin{enumerate}[(i)]
	 	\item $(\vec{u}^{\varepsilon})_{\varepsilon > 0}$ two-scale converges to $\vec{u}$. \qquad\qquad (ii) $(c^{\varepsilon})_{\varepsilon > 0}$ two-scale converges to $c$.
	 	\item[(iii)] $(w^{\varepsilon})_{\varepsilon > 0}$ two-scale converges to $w$. \qquad\;\;\;\quad (iv) $(P^{\varepsilon})_{\varepsilon > 0}$ two-scale converges to $P$.
	 	\item[(v)] $( \varepsilon  \nabla_x w^{\varepsilon})_{\varepsilon > 0}$ two-scale converges to $\nabla_{y} w$. \quad (vi) $(\varepsilon \nabla_{x} \vec{u}^{\varepsilon})_{\varepsilon > 0}$ two-scale converges to $\nabla_{y} \vec{u}$.
	 	\item[(vii)] $(\nabla_{x} c^{\varepsilon})_{\varepsilon > 0}$ two-scale converges to $\nabla_{x} c + \nabla_{y} c_1$.
	 \end{enumerate} 
\end{lemma}
\begin{proof}
	 The convergences follow from the estimates \eqref{5.1},  lemma 2.7 and lemma 2.9.
\end{proof}


In the next lemma we will discuss the convergence of nonlinear terms for $\varepsilon \rightarrow 0$.

\begin{lemma}
	The following convergence results hold:
	\begin{enumerate}[(i)]
		\item $(c^{\varepsilon})_{\varepsilon > 0}$ is strongly convergent to $c$ in $L^2(S \times \Omega)$. Thus, $\mathcal{T}^\varepsilon (c^{\varepsilon})$ converges to $c$ strongly in $L^2 (S \times \Omega \times Y)$, i.e. $(c^{\varepsilon})_{\varepsilon > 0}$ is strongly two-scale convergent to $c$.
		\item $\mathcal{T} ^{\varepsilon} \vec{u}^{\varepsilon}$ is weakly convergent to $\vec{u}$ in $L^2( S \times \Omega \times Y)^n$, i.e. $(\vec{u}^{\varepsilon})_{\varepsilon > 0}$ is weakly two-scale convergent to $\vec{u}$.
		\item $\mathcal{T} ^{\varepsilon} [ \varepsilon \nabla_{x} w^{\varepsilon}]$  converges to $ \nabla_{y} w $ weakly in $L^2(S \times \Omega \times Y)^n$, i.e. $\varepsilon \nabla_{x} w^{\varepsilon}$ is weakly two-scale convergent to $\nabla_{y} w$.
		\item The nonlinear terms $ f (c^{\varepsilon})$, $ \varepsilon c^{\varepsilon} \nabla_{x} w^{\varepsilon}$ and $   c^{\varepsilon} \vec{u}^{\varepsilon}$ two-scale converge to $ f(c) $, $ c \nabla_{y} w  $ and $  c \vec{u} $.
	\end{enumerate}
\end{lemma}
\begin{proof} 
We will prove step by step. From estimate \eqref{5.1} for $(c^{\varepsilon})_{\varepsilon > 0}$ and theorem 2.1 in \cite{MZ11}, there exists a subsequence of $(c^{\varepsilon})_{\varepsilon > 0}$, still denoted by same symbol, such that $(c^{\varepsilon})_{\varepsilon > 0}$ is strongly convergent to a limit $c$. The rest of (i) and the proofs of (ii) and (iii) follow from lemma 2.1.5.

(iv) From (i) we have the strong convergence of $c^\eps$ and by lemma 2.1.5 $T^\eps c^\eps $ is strongly convergent to $c$. Now, by corollary on page 53 in \cite{Yos70}, there exists a subsequence $c^\eps$ (denoted by same symbol) such that $T^\eps c^\eps $ is pointwise convergent to $c$, i.e. $\lim_{\eps\to 0} T^\eps c^\eps(t,x,y)=c(t,x)$ for a.e. $(t,x,y)\in S\times\Omega\times Y$. This gives that $f(T^\eps c^\eps)=[T^\eps c^\eps]^3+T^\eps c^\eps$ is pointwise convergent to $c^3+c$. We note that $||T^\eps c^\eps||_{L^4(S\times\Omega\times Y)}\le [\eps|Y|]^{\frac{1}{4}}||c^\eps||_{L^4(S\times \Omega\times Y)}\le C\; \forall \eps>0$. Again, $|f(T^\eps c^\eps)|\le |T^\eps c^\eps|^3+|T^\eps c^\eps|=:g(T^\eps c^\eps)$. Clearly, $g\ge 0$ and $g\in L^1(S\times\Omega)$. Therefore, by Lebsegue dominated convergence theorem $\lim_{\eps\to 0}\int_{S\times \Omega\times Y}f(T^\eps u^\eps)\phi \,dx,\,dy\,dt=\int_{S\times \Omega\times Y}c(t,x)\phi(t,x)\,dx\,dt=\int_{S\times \Omega}c(t,x)\phi(t,x)\,dx\,dt$, where $\phi\in L^2(S\times \Omega)$ and $|Y|=1$. Since $T^\eps f(c^\eps)=f(T^\eps c^\eps)\overset{w}{\rightharpoonup}f(c)$, by lemma 2.1.3(v) $f(c^\eps)\overset{2}{\rightharpoonup}f(c)$.

Next, we show that $ \varepsilon c^{\varepsilon} \nabla_{x} w^{\varepsilon} \overset{2}{\rightharpoonup}  c \nabla_{y} w  $, i.e. $ T^\eps(c^{\varepsilon}\varepsilon \nabla_{x} w^{\varepsilon})= \mathcal{T}^{\varepsilon} [c^{\varepsilon}] ~  \mathcal{T}^{\varepsilon} [\varepsilon \nabla_{x} w^{\varepsilon}]  \overset{w}{\rightharpoonup} c \nabla_{y} w $ as stated in lemma 2.1.3(v). We see that $\norm{  \mathcal{T}^{\varepsilon} [\varepsilon \nabla_x w^{\varepsilon}] }_{L^2(S \times \Omega \times Y )} = \eps.\varepsilon \norm{ \nabla_x w^{\varepsilon} }_{L^2(S \times \Omega)} \leq C \; \forall \eps>0$. For a $ \eta \in L^{\infty} ( S \times \Omega ; L_{\#}^{\infty}(Y))^n $ we estimate
\begin{align}
\int_{ S \times \Omega \times Y} \mathcal{T}^{\varepsilon} [c^{\varepsilon}] ~  \mathcal{T}^{\varepsilon} [\varepsilon \nabla_{x} w^{\varepsilon}] \cdot \eta \,dt \,dx \,dy  -  \int_{ S \times \Omega \times Y} c \nabla_{y} w \cdot \eta \,dt \,dx \,dy,   
\leq J_1 + J_2
\end{align}
where 
\begin{align}
 J_1  &=  \abs{\int_{ S \times \Omega \times Y} ( \mathcal{T}^{\varepsilon} [c^{\varepsilon}] - c ) ~  \mathcal{T}^{\varepsilon} [\varepsilon \nabla_{x} w^{\varepsilon}] \cdot \eta \,dt \,dx \,dy} 
\label{5.6a}
 \\
& \leq
 \norm{ \mathcal{T}^{\varepsilon} [c^{\varepsilon}] - c  }_{L^2(S \times \Omega \times Y )} \norm{ \mathcal{T}^{\varepsilon} [\varepsilon \nabla_{x} w^{\varepsilon}] }_{L^2(S \times \Omega \times Y )} \norm{ \eta }_{L^{\infty}(S \times \Omega \times Y )}
 \notag
 \\
& \rightarrow 0, \quad \text{ since $ \mathcal{T}^{\varepsilon} c^{\varepsilon}  \rightarrow c $ in $ L^2(S \times \Omega \times Y ) $ strongly.  }
\label{5.6b}
\end{align}
and
\begin{align}
J_2  &=  \abs{\int_{ S \times \Omega \times Y} ( \mathcal{T}^{\varepsilon} [\varepsilon \nabla_{x} w^{\varepsilon}]  -  \nabla_y w  )\cdot c\eta \,dt \,dx \,dy} \rightarrow 0 \quad \text{ as $ \varepsilon \rightarrow 0 $. }
\label{5.7a}
\end{align}
This implies that $T^\eps(c^{\varepsilon}\varepsilon \nabla_{x} w^{\varepsilon})$ is weakly convergent to $c \nabla_{y} w $ in $L^1(S\times\Omega\times Y).$ In other words, $\varepsilon c^{\varepsilon} \nabla_{x} w^{\varepsilon}\overset{2}{\rightharpoonup}  c \nabla_{y} w$ in $ L^1(S\times\Omega\times Y)$. 
Following the similar arguments we get $ c^{\varepsilon} \vec{u}^\varepsilon  \overset{2}{\rightharpoonup} c \vec{u} $.	
\end{proof}


\begin{proof}[Proof of Theorem 2.2]
(i) We first obtain the macroscopic description of the Cahn-Hilliard equations. We choose the function $\phi  \in C_0^{\infty} ( S \times \Omega ; C_{\#}^{\infty}(Y) )$ in \eqref{2.1b}, then
\begin{align}  
 & \int_{S} \langle \partial_t c^{\varepsilon} , \phi \rangle \,dt -  \varepsilon   \int_{S \times \Omega_p^ \varepsilon} c^{\varepsilon} \vec{u^{\varepsilon}}  \cdot \nabla    \phi \,dx \,dt +   \varepsilon^2  \int_{S \times \Omega_p^ \varepsilon} \nabla
 w^{\varepsilon} \cdot \nabla \phi \,dx \,dt = 0 \notag \\
\implies &\int_{S} \langle \partial_t c^{\varepsilon} , \phi \rangle \,dt   -    \int_{S \times \Omega_p^ \varepsilon} c^{\varepsilon} \vec{u^{\varepsilon}}  \cdot \varepsilon \nabla    \phi \,dx \,dt +  \int_{S \times \Omega_p^ \varepsilon}   \varepsilon  \nabla
w^{\varepsilon} \cdot   \varepsilon  \nabla \phi \,dx \,dt = 0\notag \\
 \implies&- \int_{S \times \Omega }   \chi(\frac{x}{\varepsilon}) c^{\varepsilon}(t, x)  \partial_t \phi (t, x, \frac{x}{\varepsilon})  \,dx \,dt 
 \notag\\
 &-  \int_{S \times \Omega } \chi(\frac{x}{\varepsilon})   c^{\varepsilon}(t, x)\vec{u^{\varepsilon}}(t, x)  \cdot    \Big(  \nabla_y \phi (t, x, y) +  \varepsilon \nabla_x \phi (t, x, y) \Big) \,dx \,dt \notag\\
&+  \int_{S \times \Omega } \chi(\frac{x}{\varepsilon})  \varepsilon  \nabla_x
w^{\varepsilon}(t, x) \cdot \Big(  \nabla_y \phi (t, x, y) +  \varepsilon \nabla_x \phi (t, x, y) \Big) \,dx \,dt  
 = 0
\end{align}
We pass $\varepsilon \rightarrow 0$ in the two-scale sense. The terms containing $\varepsilon$ are bounded and the limits converge to $0$. Hence, we get
\begin{small}
\begin{align} 
- \int_{S \times \Omega \times Y_p }  c(t, x)  \partial_t \phi (t, x, y) & \,dx \,dy \,dt -  \int_{S \times \Omega \times Y_p }   c(t, x)  \vec{u}(t, x, y)   \cdot \nabla_y \phi (t, x, y) \,dx \,dy \,dt 
\notag\\
& -  \int_{S \times \Omega \times Y_p }  \nabla_y
w(t, x, y)   \cdot \nabla_y \phi (t, x, y) \,dx \,dy \,dt  
= 0  \label{6.6}
\end{align} 
\end{small}
The strong form of \eqref{6.6} is 
\begin{align}
& \partial_t c(t, x) + \nabla_y \cdot  \{c(t, x)  \vec{u}(t, x, y) \} = \Delta_y w(t, x, y) \quad \text{  in $ S \times \Omega \times Y_p $.}
  \label{6.8}
\end{align} 

Now by choosing a function $\psi$  defined as $ \psi = \psi (t, x, \frac{x}{\varepsilon}) = \psi_0 (t, x) + \varepsilon \psi_1 (t, x, \frac{x}{\varepsilon})$ in \eqref{2.1c}, where the functions $\psi_0 \in C_0^{\infty} ( S \times \Omega )$ and $ \psi_1 \in C_0^{\infty} ( S \times \Omega ; C_{\#}^{\infty}(Y) ) $, we get 
\begin{align} 
&\int_{S} \langle w^{\varepsilon} , \psi \rangle \,dt =    \int_{S \times \Omega_p^ \varepsilon} \nabla c^{\varepsilon} \cdot \nabla \psi \,dx \,dt +  \int_{S} \langle  f(c^{\varepsilon}) , \psi \rangle \,dt 
\notag\\  
\implies & \int_{S \times \Omega_p^ \varepsilon}  w^{\varepsilon}(t, x) \psi(t, x, \frac{x}{\varepsilon})  \,dx \,dt  =    \int_{S \times \Omega_p^ \varepsilon} f(c^{\varepsilon}(t, x)) \psi(t, x, \frac{x}{\varepsilon}) \,dx \,dt
\notag\\
& +
\int_{S \times \Omega_p^ \varepsilon}  \nabla c^{\varepsilon}(t, x) \cdot \Big( \frac{1}{\varepsilon} \underbrace{\nabla_y \psi_0 (t, x)}_\text{ \mbox{=0}}  + \nabla_x \psi_0 (t, x) + \nabla_y  \psi_1 (t, x, \frac{x}{\varepsilon}) + \varepsilon  \nabla_x  \psi_1 (t, x, \frac{x}{\varepsilon})  \Big) \,dx \,dt \notag\\
 \implies &\int_{S \times \Omega }  \chi(\frac{x}{\varepsilon}) w^{\varepsilon}(t, x) \psi(t, x, \frac{x}{\varepsilon})  \,dx \,dt  =    \int_{S \times \Omega }  \chi(\frac{x}{\varepsilon}) ~ f(c^{\varepsilon}(t, x)) \psi(t, x, \frac{x}{\varepsilon}) \,dx \,dt
\notag\\
& +
 \int_{S \times \Omega } \chi(\frac{x}{\varepsilon})  \nabla c^{\varepsilon}(t, x) \cdot \Big( \nabla_x \psi_0 (t, x) + \nabla_y  \psi_1 (t, x, \frac{x}{\varepsilon}) + \varepsilon  \nabla_x  \psi_1 (t, x, \frac{x}{\varepsilon})  \Big) \,dx \,dt .
\label{6.14}
\end{align}

Next we pass $\varepsilon \rightarrow 0$ in the two-scale sense in \eqref{6.14} and get 
\begin{align}
\int_{S \times \Omega \times Y_p}    w(t, & x, y)  \psi(t, x, y)  \,dx \,dy \,dt =  
\int_{S \times \Omega \times Y_p}  f(c (t, x)) \psi(t, x, y) \,dx \,dy \,dt
\notag\\
& +  
\int_{S \times \Omega \times Y_p}  \Big( \nabla_x c (t, x)  +  \nabla_y c_1(t, x, y) \Big) \cdot   \Big( \nabla_x \psi_0 (t, x) + \nabla_y  \psi_1 (t, x, y)\Big) \,dx \,dy \,dt  .
\label{6.15}
\end{align}  

Setting $\psi_1 = 0$ in \eqref{6.15} yields 
\begin{align*}
\int_{S \times \Omega \times Y_p}    w(t, & x, y)  \psi(t, x, y)  \,dx \,dy \,dt =  
\int_{S \times \Omega \times Y_p}  f(c (t, x)) \psi(t, x, y) \,dx \,dy \,dt
\notag\\
& +  
\int_{S \times \Omega \times Y_p}  \Big( \nabla_x c (t, x)  +  \nabla_y c_1(t, x, y) \Big) \cdot    \nabla_x \psi_0 (t, x)   \,dx \,dy \,dt 
\label{6}
\end{align*} 
This gives
\begin{equation} \label{6.16}
  w(t, x, y) = -  \Delta_x c (t, x ) -   \nabla_x \cdot \nabla_y  c_1(t, x, y) +  f(c (t, x )) \quad \text{in $ S \times \Omega \times Y_p $.}
\end{equation}
Again, setting $ \psi_0 = 0 $ in \eqref{6.15} gives
\begin{align}
\int_{S \times \Omega \times Y_p}  \Big( \nabla_x c (t, x )  +  \nabla_y c_1(t, x, y) \Big) \cdot   \nabla_y  \psi_1 (t, x, y)  \,dx \,dy \,dt  = 0.
\label{6.b}
\end{align} 

One can easily deduce from \eqref{6.b}
\begin{align} 
   \nabla_y  \cdot    \nabla_y c_1(t, x, y)  = 0 \quad \text{in $ S \times \Omega \times Y_p $.}
\label{6.c} 
\end{align} 
(ii) We proceed with the homogenization of Stokes equations. We choose the test functions $ \eta \in C_0^{\infty} ( \Omega ; C_{\#)}^{\infty} (Y) ) ^n $ and $ \xi \in C_0^{\infty} ( S ) $ and proceed as in \cite{bavnas2017homogenization}. Then, using lemma (5.1) and lemma (5.2) we obtain for $ \varepsilon \rightarrow 0 $
\begin{align}
& \lim_{ \varepsilon \rightarrow 0 } \int_{ S \times \Omega_p^{\varepsilon} } P ^{\varepsilon} (t, x) \Big\{  \nabla_x   \cdot \eta ( x, \frac{x}{\varepsilon} )   +  \frac{1}{\varepsilon} \nabla_y \cdot \eta ( x, \frac{x}{\varepsilon} ) \Big\}  \partial_t \xi (t) \,dx \,dy \,dt0 
\notag\\
&= \int_{ S \times \Omega \times Y_p} P(t,x,y) \nabla_y \cdot \eta ( x, y ) \partial_t \xi (t) \,dx \,dy \,dt 
\notag\\
&= 0
\label{6.17}
\end{align}
From \eqref{6.17} we deduce that the two-scale limit of the pressure $P$ is independent of $y$, i.e., $P(t, x) \in L_0^2 ( S \times \Omega ) $. Next, we take the function $ \eta \in C_0^{\infty} ( \Omega ; C_{\#}^{\infty} (Y) ) ^n $ such that $ \nabla_y \cdot \eta ( x, y ) = 0 $, then
\begin{align}
\mu  \varepsilon^2 \int_{ S \times \Omega_p^{\varepsilon} } \nabla  \vec{u}^{\varepsilon} (t, x) : \nabla \eta ( x, y )  \xi (t) &  \,dx  \,dt
+
\int_{ S \times \Omega_p^{\varepsilon} }  P^{\varepsilon}(t, x) \nabla \cdot  \eta ( x, y ) \partial_t \xi (t)  \,dx  \,dt 
\notag \\
& =
- \lambda \varepsilon \int_{S \times \Omega_p^{\varepsilon} } c^{\varepsilon} ( t, x ) \nabla w^{\varepsilon} ( t, x ) \cdot \eta ( x, y )  \xi (t)  \,dx  \,dt  
\label{6.18}. 
\end{align}
Using the extensions of solution to $\Omega$ with the same notation, we get
\begin{align}
\mu  \varepsilon^2 \int_{ S \times \Omega } \chi ( \frac{x}{ \varepsilon } )  \nabla  \vec{u}^{\varepsilon} (t, x) : \nabla \eta ( x, y ) & ~ \xi (t)  \,dx   \,dt 
+
\int_{ S \times \Omega }  \chi ( \frac{x}{ \varepsilon } ) P^{\varepsilon}(t, x) \nabla \cdot  \eta ( x, y ) ~ \partial_t \xi (t) \,dx  \,dt
\notag \\
& =
- \lambda ~ \varepsilon \int_{ S \times \Omega } \chi ( \frac{x}{ \varepsilon } ) c^{\varepsilon} ( t, x ) \nabla w^{\varepsilon} ( t, x ) \cdot \eta ( x, y ) ~ \xi (t)   \,dx  \,dt
\notag, 
\\
\implies \mu   \int_{ S \times \Omega } \chi ( \frac{x}{ \varepsilon } ) ~ \varepsilon \nabla  \vec{u}^{\varepsilon} (t, x) : \varepsilon \nabla \eta ( x, y ) & ~ \xi (t)  \,dx   \,dt 
+
\int_{ S \times \Omega }  \chi ( \frac{x}{ \varepsilon } ) P^{\varepsilon}(t, x) \nabla \cdot  \eta ( x, y ) ~ \partial_t \xi (t) \,dx  \,dt
\notag \\
& =
- \lambda   \int_{ S \times \Omega } \chi ( \frac{x}{ \varepsilon } ) ~ c^{\varepsilon} ( t, x ) ~ \varepsilon \nabla w^{\varepsilon} ( t, x ) \cdot \eta ( x, y ) ~  \xi (t)   \,dx  \,dt
\notag, 
\end{align}
\begin{align}
\implies \mu   \int_{ S \times \Omega } \chi ( \frac{x}{ \varepsilon } ) ~ \varepsilon \nabla  \vec{u}^{\varepsilon} (t, x) & :  \{ \nabla_y \eta  ( x, y )  + \varepsilon \nabla_x \eta  ( x, y )  \}   \xi (t)  \,dx   \,dt 
\notag \\
& 
+ \int_{ S \times \Omega }  \chi ( \frac{x}{ \varepsilon } ) P^{\varepsilon}(t, x) \Big\{ \frac{1}{\varepsilon} \nabla_y \cdot \eta  ( x, y )  +  \nabla_x  \cdot \eta  ( x, y )  \Big\}  ~ \partial_t \xi (t) \,dx  \,dt
\notag \\
&=
- \lambda   \int_{ S \times \Omega } \chi ( \frac{x}{ \varepsilon } ) ~ c^{\varepsilon} ( t, x ) ~ \varepsilon \nabla w^{\varepsilon} ( t, x ) \cdot \eta ( x, y ) ~  \xi (t)   \,dx  \,dt
\notag,
\end{align}
\begin{align}
\implies \mu   \int_{ S \times \Omega } \chi ( \frac{x}{ \varepsilon } ) ~ \varepsilon \nabla  \vec{u}^{\varepsilon} (t, x) & :  \{ \nabla_y \eta  ( x, y )  + \varepsilon \nabla_x \eta  ( x, y )  \}   \xi (t)  \,dx   \,dt 
\notag \\
& 
+ \int_{ S \times \Omega }  \chi ( \frac{x}{ \varepsilon } ) P^{\varepsilon}(t, x) \{ \nabla_x  \cdot \eta  ( x, y )  \}    ~ \partial_t \xi (t) \,dx  \,dt
\notag \\
&=
- \lambda   \int_{ S \times \Omega } \chi ( \frac{x}{ \varepsilon } ) ~ c^{\varepsilon} ( t, x ) ~ \varepsilon \nabla w^{\varepsilon} ( t, x ) \cdot \eta ( x, y ) ~  \xi (t)   \,dx  \,dt
\label{6.22},
\end{align}

Next, we pass $ \varepsilon \rightarrow 0 $ in the two-scale sense. The terms containing $\varepsilon$ are bounded and the limits converge to 0. Hence, we get
\begin{align}
\implies \mu   \int_{ S \times \Omega \times Y_p} \nabla_y \vec{u} (t, x, y) : \nabla_{y} \eta ( x, & y )  \xi (t)  \,dx \,dy \,dt 
 +
 \int_{ S \times \Omega \times Y_p} P(t, x) \nabla_x \cdot \eta ( x, y ) \partial_t \xi (t) \,dx  \,dy \,dt 
\notag \\
 & =
- \lambda \int_{ S \times \Omega \times Y_p} c ( t, x ) \nabla_y w ( t, x, y ) \cdot \eta ( x, y )  \xi (t)  \,dx \,dy \,dt  
\label{6.23}. 
\end{align} 

We follow the existence of a pressure $P_1 \in L^{\infty} ( S ; L^2_0 ( \Omega ; L^2_{\#} ( Y_p ) ) ) $ and two-scale convergence results as in \cite{bavnas2017homogenization} for the last step of our proof.
\begin{small}
\begin{align}
\implies  \mu   \int_{ S \times \Omega \times Y_p} \nabla_y \vec{u} (t, x, y) : &\nabla_{y} \eta ( x,  y )  \xi (t)  \,dx \,dy \,dt 
+
\int_{ S \times \Omega \times Y_p} P(t, x) \nabla_x \cdot \eta ( x, y ) \partial_t \xi (t) \,dx  \,dy \,dt 
\notag \\
& +
\int_{ S \times \Omega \times Y_p} P_1(t, x, y) \nabla_y \cdot \eta ( x, y ) \partial_t \xi (t) \,dx  \,dy \,dt  
\notag \\ 
& =
- \lambda \int_{ S \times \Omega \times Y_p} c ( t, x  ) \nabla_y w ( t, x, y ) \cdot \eta ( x, y )  \xi (t)  \,dx \,dy \,dt  
\label{6.24}. 
\end{align} 
\end{small}
for all  $ \eta \in C_0^{\infty} ( \Omega ; C_{\#}^{\infty} (Y) ) ^n $ and $  \xi \in C_0^{\infty} ( S ) $.

From \eqref{6.24}, we obtain 
\begin{equation}
- \mu \Delta_y \vec{u} ( x, y) + \nabla_x p ( x ) + \nabla_y p_1 (x, y)
=
-  \lambda ~  c ( x  ) ~ \nabla_y w ( x, y ) \quad \text{in $ S \times \Omega \times Y_p $.}
\end{equation}

\end{proof}

\section{Conclusion}
In this paper, we considered a mixture of two fluids in a heterogeneous porous medium, where the fluids in the pore space were separated by an interface of thickness of $\lambda$. The modeling of such processes lead to a strongly coupled system of Stokes-Cahn-Hilliard equations, i.e. a system of parabolic and elliptic equations. The model equation are taken into account at the microscopic scale. We chose the case where $\alpha =2, \beta =1$ and $\gamma =0$ in model equations (1.5). The effect of surface tension is incorporated in the model and the interfacial layer (of thickness $\lambda$) is assumed to be $\eps$ dependent. We first showed the well-posedness of the model at the micro scale and via energy method we also obtained several a-priori estimates. The existence of solution follows mainly from Galerkin's approximation and can be found in \cite{Fen06,SL03,FHC07}. The estimates together with two-scale convergence and periodic unfolding were used to upscale the model equation from micro to macro scale. The model and analysis proposed in this paper can be further utilized to generalize the model to a three or more multicomponent fluid mixtures. This question will be addressed elsewhere. 
\section{Acknowledgment}
The first author would like to thank IIT Kharagpur for providing the funding for her PhD position. 
\small 
\bibliographystyle{plain}
	\bibliography{paper1}

\end{document}